\documentclass{amsart}
\usepackage{graphicx}
\usepackage{amssymb,amscd,amsthm,amsxtra}
\usepackage{latexsym}
\usepackage{epsfig}
\usepackage{color}
\newtheorem{thm}{Theorem}[section]
\newtheorem{cor}[thm]{Corollary}
\newtheorem{lem}[thm]{Lemma}
\newtheorem{prop}[thm]{Proposition}
\theoremstyle{definition}
\newtheorem{defn}[thm]{Definition}
\theoremstyle{remark}
\newtheorem{rem}[thm]{Remark}
\numberwithin{equation}{section}
\newcommand{\R}{\mathbb R}
\newcommand{\be}{\begin{equation}}
\newcommand{\ee}{\end{equation}}

\newcommand{\ep}{\eps}

\newcommand{\eps}{\varepsilon}

\newcommand{\comment}[1]{}
\begin{document}

\title[Regularity of flat free boundaries for the $p(x)$-Laplacian]{Regularity of flat free boundaries for a $p(x)$-Laplacian problem with right hand side}

\author[Fausto Ferrari]{Fausto Ferrari}
\address{Dipartimento di Matematica dell'Universit\`a di Bologna, Piazza di Porta S. Donato, 5, 40126 Bologna, Italy.}
\email{\tt fausto.ferrari@unibo.it}
\author[Claudia Lederman]{Claudia Lederman}
\address{IMAS - CONICET and Departamento  de
Ma\-te\-m\'a\-ti\-ca, Facultad de Ciencias Exactas y Naturales,
Universidad de Buenos Aires, (1428) Buenos Aires, Argentina.}
\email{\tt  clederma@dm.uba.ar}
\thanks{F. F. was partially supported by INDAM-GNAMPA 2019 project: {\it Proprietà di regolarità delle soluzioni viscose con applicazioni a problemi di frontiera libera.} }
\thanks{C. L. was partially supported by the project GHAIA Horizon 2020 MCSA RISE 2017 programme grant 777822 and  by the grants CONICET PIP 11220150100032CO 2016-2019,  UBACYT 20020150100154BA and ANPCyT PICT 2016-1022. C. L. wishes to thank the Department of Mathematics of the University of Bologna, Italy, for the kind hospitality.}
%35R35 Free boundary problems for PDE
%35B65 Smoothness and regularity of solutions of PDE
%35J60 Nonlinear elliptic equations
%35J70 Degenerate elliptic equations

\begin{abstract}
We consider viscosity solutions to a one-phase free boundary problem for the $p(x)$-Laplacian with non-zero right hand side. We apply the tools developed in \cite{D} to prove that flat free boundaries are $C^{1,\alpha}$.
Moreover, we obtain some new results for the operator under consideration that are of independent interest.
\end{abstract}

% ----------------------------------------------------------------
\maketitle
% ----------------------------------------------------------------

\section{Introduction and main results}\label{section1}

In this paper we study a one-phase  free boundary problem governed by the $p(x)$-Laplacian  with non-zero right hand side. More precisely, we denote  by
$$
\Delta_{p(x)}u:=\mbox{div} (|\nabla u|^{p(x)-2}\nabla u),
$$
where $p$ is a  function such that $1<p(x)<+\infty$. Then our problem is the following:
\begin{equation}  \label{fb}
\left\{
\begin{array}{ll}
\Delta_{p(x)} u = f, & \hbox{in $\Omega^+(u):= \{x \in \Omega : u(x)>0\}$}, \\
\  &  \\
|\nabla u|= g, & \hbox{on $F(u):= \partial \Omega^+(u) \cap
\Omega.$} 
\end{array}
\right.
\end{equation}
Here $\Omega  \subset \mathbb{R}^n$ is a bounded domain,  $p\in C^1(\Omega)$,  
$f\in C(\Omega)\cap L^{\infty}(\Omega)$ and $g\in C^{0, \beta}(\Omega),$ $g\geq 0.$ 

This problem comes out  naturally from limits of a singular perturbation problem with
forcing term as in \cite{LW1}, where the authors analyze solutions to  
\eqref{fb}, arising in the study of flame propagation with nonlocal and electromagnetic 
effects. On the other hand, \eqref{fb} appears by minimizing the following functional
\begin{equation}\label{AC-energy}
\mathcal{E}(v)=\int_{\Omega}\left(\frac{|\nabla v|^{p(x)}}{p(x)}+Q^2(x)\chi_{\{v>0\}}+f(x)v\right)dx
\end{equation}
studied in \cite{LW3}, as well as in the seminal paper by Alt and Caffarelli \cite{AC} in the case $p(x)\equiv 2$ and $f\equiv 0.$ We refer also to \cite{LW4}, where 
\eqref{fb} appears in the study of an optimal design problem.

We are interested in the regularity of the free boundary for viscosity solutions of (\ref{fb}). This problem has been already faced in \cite{LW2} for weak solutions with the aid of the techniques developed  in \cite{AC}. 

In the present work we are following the strategy introduced in the important paper by De Silva \cite{D}, that was inspired by \cite{S}, for one-phase problems and linear non-divergence operators.  \cite{D} was further extended to two-phase problems in different settings, see 
\cite{DFS1, DFS2, DFS3}.   The same technique was applied to the $p$-Laplace operator ($p(x)\equiv p$ in \eqref{fb}) for the one phase case, with $p\ge 2$, in \cite{LR}. See also \cite{LT}.

In the linear homogeneous case, $f\equiv 0,$ \eqref{fb} was studied for viscosity solutions in the pioneer works by
Caffarelli \cite{C1,C2}. 
The results in \cite{C1,C2} have been widely generalized to different
classes of homogeneous elliptic problems. See for example \cite{CFS, FS1,
FS2} for linear operators, \cite{AF, F1, F2, Fe1, W1, W2, RT} for fully nonlinear
operators and \cite{LN1, LN2} for the $p$-Laplacian. See also \cite{ART}.

	As already mentioned, problem \eqref{fb} was originally studied in the linear homogeneous case in \cite{AC}, associated to \eqref{AC-energy}. These techniques were generalized to the linear case with $f\not\equiv 0$ in \cite{GS, Le}. In the homogeneous case,   to a quasilinear uniformly elliptic situation \cite{ACF}, to the $p$-Laplacian \cite{DP}, to an Orlicz setting \cite{MW} and to the $p(x)$-Laplacian with $p(x)\ge 2$ \cite{FMW}. Finally,  \eqref{fb} with $1<p(x)<\infty$ and $f\not\equiv 0$ was dealt with in \cite{LW2}.

\medskip

In this paper we show that flat free boundaries of viscosity solutions to
\eqref{fb} are $C^{1,\alpha}.$ In the forthcoming work \cite{FL} we  prove that Lipschitz free boundaries of viscosity solutions to
\eqref{fb} are $C^{1,\alpha}.$

\smallskip
 
 Our main result is the following (for the precise definition of viscosity solution to \eqref{fb} we refer to Section \ref{section2})
\begin{thm}[Flatness implies $C^{1,\protect\alpha}$]
\label{flatmain1} Let
$u$ be a viscosity solution to \eqref{fb}
in $B_1$. Assume that  $0\in F(u),$ $g(0)=1$ and $p(0)=p_0.$   
There exists a universal constant $\bar{\varepsilon}>0$ such that, if the graph of $u$ is $\bar{\varepsilon}-$flat in $B_1,$  in the direction $e_n,$ that is
\begin{equation}  \label{cflat}
(x_n-\bar{\varepsilon})^+\leq u(x)\leq (x_n+\bar{\varepsilon})^+, \quad x\in B_1,
\end{equation}
and
\begin{equation}  \label{pflat}
\|\nabla p\|_{L^{\infty}(B_1)}\leq \bar{\varepsilon},\quad \|f\|_{L^{\infty}(B_1)}\leq \bar{\varepsilon}, \quad [g]_{C^{0,\beta}(B_1)}\leq \bar{\varepsilon},
\end{equation}
 then $F(u)$ is $C^{1,\alpha}$ in $
B_{1/2}$.
\end{thm}

\medskip

In addition to the assumptions already stated above, we suppose that
\begin{equation}\label{p-lip} \nabla p \in L^{\infty}(\Omega)
 \end{equation}
and that there exist positive numbers  $p_{\min},p_{\max},$ such that 
\begin{equation}\label{exponentsize}1<p_{\min}\leq p(x)\leq p_{\max}<\infty.
 \end{equation}
In Theorem \ref{flatmain1} the constants $\bar{\varepsilon}$ and $\alpha$ depend only on $p_{\min}$, $p_{\max}$  and $n$ (the dimension of the space).

The proof of Theorem \ref{flatmain1} is based on an
improvement of flatness, obtained via a compactness argument which linearizes the
problem into a limiting one.  The key tool is a geometric Harnack inequality
that localizes the free boundary well, and allows the rigorous passage to
the limit.

Let us point out that carrying out, for the inhomogeneous $p(x)$-Laplace operator, the strategy devised in \cite{D}   required the development of new tools. In fact, the $p(x)$-Laplacian  is a nonlinear operator that appears naturally in divergence form from minimization problems, i.e., in the form ${\rm div}A(x,\nabla u)=f(x)$, with
\begin{equation*}\label{Fan-1.6}
\lambda |\eta|^{p(x)-2}|\xi|^2\le\sum_{i,j=1^n}\frac{\partial A_i}{\partial \eta_j}(x,\eta)\xi_i\xi_j\le \Lambda |\eta|^{p(x)-2}|\xi|^2, \quad \xi\in\R^n.
\end{equation*}
This operator is singular in the regions where $1<p(x)<2$  and  degenerate in the ones where $p(x)>2$. 

Some  results for this type of operators we needed to use to achieve our goals are available in the literature for weak solutions (in the sense of Definition \ref{defnweak}  in Section \ref{section3}). These results are Harnack inequality (see \cite{Wo}) and $C^{1,\alpha}$ estimates (see \cite{Fan} and \cite{FanZ}).
However, the program followed in \cite{D} relies on solutions of the corresponding equations in a viscosity sense (see \cite{CIL}).

The equivalence between weak and viscosity solutions of $\Delta_{p(x)}u=f$ was proved in \cite{JJ,JLM,  MO} in the case of the $p$-Laplacian (i.e., for $p(x)\equiv p$) and in \cite{JLP} in the case of the homogeneous $p(x)$-Laplacian (i.e., for $f\equiv 0$). To our knowledge there is no such result in the literature for the inhomogeneous $p(x)$-Laplacian.

Hence, in order to proceed with the arguments in \cite{D}, we prove in Theorem \ref{weak-is-visc} that weak solutions of $\Delta_{p(x)}u=f$  are indeed viscosity solutions. This new result is of independent interest, since it may be applied in other contexts.

On the other hand, the approach in \cite{D} requires the use of  barriers of the type 
$w(x)=c_1|x-x_0|^{-\gamma}-c_2$, together with suitable  modification of them. In the present work we are able to employ the same kind of barriers. Showing that they are also appropriate  to deal with the inhomogeneous $p(x)$-Laplace operator was a nontrivial and delicate task, that we perform in Lemma \ref{barry}. Again, the difficulty relies on the nonlinear singular/degenerate nature and $x$ dependence  of our equation and also on the presence of the logarithmic term appearing in the nondivergence form of the operator (see \eqref{non-diverg}).

The results in Lemma \ref{barry} are new even for $p(x)\equiv p$ in the range $1<p<2$. These barriers, which are novel in the $p(x)$-Laplace context, are different from the ones used in the literature for this operator (see, for instance, \cite{FMW, Wo, LW4}). Consequently, our results in Lemma \ref{barry} have possible applications to other situations.

We would like to stress at this stage that partial differential equations with non-standard growth have been receiving a lot of attention and that the $p(x)$-Laplacian is a model case in this class. A list of applications of this type of operators includes
the modelling of non-Newtonian fluids,  for instance, electrorheological \cite{R} or thermorheological fluids \cite{AR}. Also non-linear elasticity \cite{Z1},  image reconstruction \cite{AMS,CLR} and the modelling of electric conductors \cite{Z2}, to cite a few.

The fact that solutions to the inhomogeneous $p(x)$-Laplacian are locally of class $C^{1,\alpha}$ plays a critical
role in the analysis of this paper. A comprehensive account for
sharp conditions for regularity of solutions of some elliptic equations  with
non-standard growth can be found in  \cite{AM} and \cite{Fan}.

We finally remark that our main result, Theorem \ref{flatmain1}, is applied in the companion paper \cite{FL} to  prove that Lipschitz free boundaries of viscosity solutions of \eqref{fb} are $C^{1,\alpha}$.

Our work is organized a follows. In Section \ref{section2} we provide notation and basic definitions, and we also present an auxiliary result on a Neumann problem which will be used in the proof of Theorem \ref{flatmain1}. In Section 3 we discuss the relationship between the different notions of solutions to $\Delta_{p(x)}u=f$ we are using. In particular, we prove Theorem \ref{weak-is-visc} which shows that weak solutions to $\Delta_{p(x)}u=f$ are viscosity solutions of the same equation. In Section 4 we prove some auxiliary results, which include Lemma \ref{barry}, concerning the existence of barrier functions for $\Delta_{p(x)}u=f$. Next, in Section 5 we prove a geometric Harnack inequality for problem \eqref{fb}. In Section \ref{sect-improv} we prove an improvement of flatness lemma. Finally, in Section \ref{section7} we prove our main result, Theorem \ref{flatmain1}. For the sake of completeness, we also include an Appendix at the end of the paper where we introduce the Sobolev spaces with variable exponent, which are the appropriate spaces to work with weak solutions of the $p(x)$-Laplacian.

\section{Basic definitions, notation and preliminaries}\label{section2}

In this section, we provide notation and basic definitions we will use throughout  our work. We also present an auxiliary result on a Neumann problem that will be applied in the paper.

\smallskip

\noindent {\it Notation.} For any continuous function $u:\Omega\subset \mathbb{R}^n\to \mathbb{R}$ we denote
\begin{equation*}
\Omega^+(u):= \{x \in \Omega : u(x)>0\},\qquad F(u):= \partial \Omega^+(u) \cap \Omega. 
\end{equation*} 
We refer to the set $F(u)$ as the {\it free boundary} of $u$, while $\Omega^+(u)$ is its {\it positive phase} (or {\it side}).

\smallskip

 Below we give the definition of viscosity solution to problem (\ref{fb}) and we deduce some consequences. In particular, we refer to the usual $C$-viscosity definition of sub/supersolution and solution of an elliptic PDE, see e.g. \cite{CIL}.

First we need the following standard notion.

\begin{defn}Given $u, \varphi \in C(\Omega)$, we say that $\varphi$
touches $u$ from below (resp. above) at $x_0 \in \Omega$ if $u(x_0)=
\varphi(x_0),$ and
$$u(x) \geq \varphi(x) \quad (\text{resp. $u(x) \leq
\varphi(x)$}) \quad \text{in a neighborhood $O$ of $x_0$.}$$ If
this inequality is strict in $O \setminus \{x_0\}$, we say that
$\varphi$ touches $u$ strictly from below (resp. above).
\end{defn}

\begin{defn}\label{defnhsol1} Let $u$ be a continuous nonnegative function in
$\Omega$. We say that $u$ is a viscosity solution to (\ref{fb}) in
$\Omega$, if the following conditions are satisfied:
\begin{enumerate}
\item $ \Delta_{p(x)} u = f$ in $\Omega^+(u)$ 
in the weak sense of Definition \ref{defnweak}, see Section \ref{section3}.
\item For every $\varphi \in C(\Omega)$, $\varphi \in C^2(\overline{\Omega^+(\varphi)})$. If $\varphi^+$ touches $u$ from below (resp.  above) at $x_0 \in F(u)$ and $\nabla \varphi(x_0)\not=0$, then 
$$|\nabla \varphi(x_0)| \leq g(x_0)\quad (\text{resp. $ \geq g(x_0)$)}.$$
\end{enumerate}
\end{defn}

Next theorem follows as a consequence of Theorem \ref{weak-is-visc} in Section \ref{section3}.

\begin{thm}\label{defnhsol2} Let $u$ be a viscosity solution to (\ref{fb}) in
$\Omega.$  Then the following conditions are satisfied:
\begin{enumerate}
\item $ \Delta_{p(x)} u = f$ in $\Omega^+(u)$ in the
viscosity sense, 
that is:
\begin{itemize}
\item[(ia)] for every $\varphi\in C^2(\Omega^+(u))$ and for every $x_0\in \Omega^+(u),$ if $\varphi$ touches $u$ from above at $x_0$ and $\nabla \varphi(x_0)\not=0,$ then $\Delta_{p(x_0)}\varphi(x_0)\geq f(x_0),$ that is, $u$ is a viscosity subsolution;
\item[(ib)] for every $\varphi\in C^2(\Omega^+(u))$ and for every $x_0\in \Omega^+(u),$ if $\varphi$ touches $u$ from below at $x_0$ and $\nabla \varphi(x_0)\not=0,$ then $\Delta_{p(x_0)}\varphi(x_0)\leq f(x_0),$ that is, $u$ is a viscosity supersolution.
\end{itemize}
\item For every $\varphi \in C(\Omega)$, $\varphi \in C^2(\overline{\Omega^+(\varphi)})$. If $\varphi^+$ touches $u$ from below (resp.  above) at $x_0 \in F(u)$ and $\nabla \varphi(x_0)\not=0$, then 
$$|\nabla \varphi(x_0)| \leq g(x_0)\quad \mbox{(resp.} \geq g(x_0)\mbox{)}.$$
\end{enumerate}
\end{thm}

\begin{rem}
If $p(x)\equiv p$ or $f\equiv 0$, then any function satisfying the conditions of Theorem \ref{defnhsol2} is a solution to \eqref{fb} in the sense of Definition \ref{defnhsol1} (see Remark \ref{equiv-not}).
\end{rem}

We introduce also the notion of comparison sub/supersolution.

\begin{defn}\label{defsub}
We say that $v \in C(\Omega)$ is a strict (comparison) subsolution (resp.
supersolution) to (\ref{fb}) in $\Omega$ if $v \in C^2(\overline{\Omega^+(v) })$,  $\nabla v\not=0$ in $\overline{\Omega^+(v) }$ and the
following conditions are satisfied:
\begin{enumerate}
\item $ \Delta_{p(x)} v  > f $ (resp. $< f $) in $\Omega^+(v)$; 
\item If $x_0 \in F(v)$, then $$|\nabla v(x_0)|>g(x_0)\quad (\text{resp. $|\nabla v(x_0)| <g(x_0)$}). $$
\end{enumerate}
\end{defn}

Notice that by the implicit function theorem, according to our definition, the free boundary of a comparison sub/supersolution is $C^2$.

\smallskip

As a consequence of the previous discussion we have

\begin{lem}\label{compare} Let $u$ be a viscosity solution to \eqref{fb} in $\Omega$. If $v$ is a strict  (comparison) subsolution to \eqref{fb} in $\Omega$ and  $u\ge v^+$ in $\Omega$ then $u>v$ in $\Omega^+(v)\cup F(v)$. Analogously, 
if  $v$ is a strict  (comparison) supersolution to \eqref{fb} in $\Omega$ and $v\ge u$ in $\Omega$ then $v>u$ in $\Omega^+(u)\cup F(u)$.
\end{lem}

\noindent {\it Notation.} From now on $B_{\rho}(x_0)\subset {\R}^n$ will denote the open ball of radius $\rho$ centered at $x_0$, and 
$B_{\rho}=B_{\rho}(0)$. A positive constant depending only on the dimension $n$, $p_{\min}$, $p_{\max}$ will be called a universal constant. We will use $c$, $c_i$ to denote small universal constants and $C$, $C_i$ to denote large universal constants.

\bigskip

The rest of the section is devoted to the study of the linearized problem associated with  our free boundary problem \eqref{fb}. That is,
the classical Neumann problem for a constant coefficient linear operator. Precisely, we consider  the following boundary
value problem:
\begin{equation}\label{Neumann_p}
  \begin{cases}
    {\mathcal L}_{p_0} \tilde u=0 & \text{in $B_\rho \cap \{x_n > 0\}$}, \\
\tilde u_n=0 & \text{on $B_\rho \cap \{x_n =0\}$}.
  \end{cases}\end{equation}
Here $1<p_{\min}\le p_0\le p_{\max}<\infty$, $\tilde u_n$ denotes the derivative in the $e_n$ direction of $\tilde u$ and
\begin{equation}\label{Lp0}
{\mathcal L}_{p_0}  u := \Delta u +(p_0-2) {\partial}_{nn} u.
\end{equation}

Theorem \ref{flatmain1} will follow via a compactness argument combined with regularity properties of solutions to \eqref{Neumann_p}, namely Theorem \ref{linearreg}.

We use the notion of viscosity solution to \eqref{Neumann_p}. We recall standard notions and a regularity result for viscosity solutions
to \eqref{Neumann_p}.

\begin{defn} \label{def-neumann}
Let $\tilde u$ be a continuous function on $B_{\rho}\cap\{x_n\ge 0\}$. We say that $\tilde u$ is a viscosity solution to \eqref{Neumann_p} if given a quadratic polynomial $P(x)$ touching $\tilde u$ from below (resp. above) at $\bar{x}\in B_\rho \cap \{x_n \ge 0\}$,
\begin{itemize}
\item[(i)] if $\bar{x}\in B_\rho \cap \{x_n > 0\}$ then ${\mathcal L}_{p_0}P\le 0$ (resp. ${\mathcal L}_{p_0}P\ge 0$), i.e.  
${\mathcal L}_{p_0}\tilde u= 0$ in the viscosity sense in $B_\rho \cap \{x_n > 0\}$;

\item[(ii)] if $\bar{x}\in B_\rho \cap \{x_n = 0\}$ then $P_n(\bar{x})\le 0$ (resp. $P_n(\bar{x})\ge 0$).
\end{itemize}
\end{defn}

\begin{rem} \label{polyn-strict} Notice that in the definition above we can choose polynomials $P$ that touch $\tilde u$ strictly from above/below. Also, it suffices to verify that (ii) holds for polynomials $\widetilde P$ with ${\mathcal L}_{p_0}\widetilde P >0$ (see \cite{D}).
\end{rem}

We will use the following regularity result for viscosity solutions to the linearized problem \eqref{Neumann_p}. For the proof we refer to Theorem 7.4 in \cite{MS}.

\begin{thm}\label{linearreg}Let $\tilde u$ be a viscosity solution to \eqref{Neumann_p} in $B_{1/2}\cap \{x_n \ge 0\}$. Then, 
$\tilde u\in C^2(B_{1/2}\cap\{x_n\ge 0\})$ and it is a classical solution to \eqref{Neumann_p}. 

Moreover, if $\|\tilde u\|_\infty \leq 1$, then there exists a constant $\bar C>0$, depending only on $n, p_{\min}$ and $p_{\max}$, such that
\be\label{lr}|\tilde u(x) - \tilde u(0) -\nabla\tilde u(0)\cdot x| \leq \bar C r^2 \quad \text{in $B_r\cap \{x_n \ge 0\}$},\ee for all $r \leq 1/4$.
\end{thm}

\section{Different notions of solutions to ${p(x)}$-Laplacian} \label{section3}

In this section we discuss the relationship between the different notions of solutions to $\Delta_{p(x)}u=f$ we are using, namely weak and viscosity solutions.

We start by observing that 
direct calculations show that, for $C^2$ functions $u$ such that $\nabla u(x)\not=0$,
\begin{equation}\label{non-diverg}
\begin{split}
&\Delta_{p(x)}u=\mbox{div} (|\nabla u|^{p(x)-2}\nabla u)\\
&=|\nabla u(x)|^{p(x)-2}\left(\Delta u+(p(x)-2)\Delta_\infty^Nu+\langle \nabla p(x),\nabla u(x)\rangle\log |\nabla u(x)| \right),
\end{split}
\end{equation}
where
$$
\Delta_\infty^Nu:=\Big\langle D^2 u(x)\frac{\nabla u(x)}{|\nabla u(x)|}\,,\,\frac{\nabla u(x)}{|\nabla u(x)|}\Big\rangle
$$
denotes the normalized $\infty$-Laplace operator.

First we need (see the Appendix for the definition of Sobolev spaces with variable exponent)
\begin{defn}\label{defnweak} 
Assume that $1<p_{\min}\le p(x)\le p_{\max}<\infty$
with  $p(x)$ Lipschitz continuous in $\Omega$ and $\|\nabla p\|_{L^{\infty}}\leq L$, for some $L>0$ and  $f\in L^{\infty}(\Omega)$.

We say that $u$
is a weak solution to $\Delta_{p(x)}u=f$ in $\Omega$ if $u\in W^{1,p(\cdot)}(\Omega)$ and,  for every  $\varphi \in
C_0^{\infty}(\Omega)$, there holds that
$$
-\int_{\Omega} |\nabla u(x)|^{p(x)-2}\nabla u \cdot \nabla
\varphi\, dx =\int_{\Omega} \varphi\, f(x)\, dx.
$$
\end{defn}

We next prove
\begin{thm} \label{weak-is-visc} Let $p$ and $f$ be as in Definition \ref{defnweak}. Assume moreover that $f\in C(\Omega)$ and $p\in C^1(\Omega)$. 

 Let $u\in W^{1,p(\cdot)}(\Omega)\cap C(\Omega)$ be a weak solution to $\Delta_{p(x)}u=f$ in $\Omega$. Then $u$ is a viscosity solution to $\Delta_{p(x)}u=f$ in $\Omega$.
\end{thm}
\begin{proof} 

Let us show that $u$ is a  viscosity supersolution to $\Delta_{p(x)}u=f$ in $\Omega$. 

{\it Step I.} We will first prove the result under the extra assumption that $f\in W^{1,\infty}(\Omega)$ and $p\in C^{1,\beta}(\Omega)$, for some $0<\beta<1$.

In fact, let $v\in C^2(\Omega)$ such that  $v$ touches $u$ from below at $x_0\in \Omega$, with $\nabla v(x_0)\not=0$. We will show that 
\begin{equation}\label{u-supers}
\Delta_{p(x_0)}v(x_0)\leq f(x_0).
\end{equation}

Let us fix $r>0$ such that $\overline{B_r(x_0)}\subset\Omega$. From Theorem 1.1 in \cite{Fan} we know that $u\in C^{1,\alpha}$ in  $\overline{B_r(x_0)}$, for some $0<\alpha<1$. We can assume that $\alpha\le \beta$.

Since $v$ touches $u$ from below at $x_0$, we know that $\nabla u(x_0)=\nabla v(x_0)\not=0$. Then, we can choose $r$ small enough so that 
$$c_1\le |\nabla u(x)| \le C_1 \quad \text{ in } B_r(x_0), \quad (c_1, C_1 \text{ positive constants}).$$

Now, arguing as in Theorem 3.2 in \cite{CL} we deduce that $u\in W^{2,2}_{\rm loc}(B_r(x_0))$ and it is a solution to the linear uniformly elliptic equation
\[
\Delta_{p(x)}u=\sum_{i,j=1}^n a_{ij}(x){u}_{x_ix_j}+\sum_{i=1}^n b_i(x){u}_{x_i}=f \quad \text{ in } B_r(x_0)
\]
where
$$a_{ij}(x)=|\nabla
u|^{p(x)-2}\Big(\delta_{ij}+(p(x)-2)\frac{{u}_{x_i}{u}_{x_j}}{|\nabla
u|^2}\Big),$$
and
$$b_i(x)=|\nabla
u|^{p(x)-2}\Big({p}_{x_i}(x)\log|\nabla u|\Big),$$ 
with
\begin{equation*}
\beta_1|\xi|^2\le \sum_{i,j=1}^n a_{ij}(x)\xi_i\xi_j\le
\beta_2|\xi|^2, \quad\forall \xi\in\R^N, \ \forall x\in B_r(x_0),
\end{equation*}
for $\beta_1, \beta_2$ positive constants. It follows (see, for instance, Theorem 9.19 in \cite{GT}) that $u\in C^{2,\alpha}$ in  $B_r(x_0)$.

Since $v$ touches $u$ from below at $x_0$, we have $\nabla u(x_0)=\nabla v(x_0)$ and $D^2u(x_0)\ge D^2v(x_0)$ and then,
\begin{equation*}
\begin{aligned}
f(x_0)&=\Delta_{p(x_0)}u(x_0)\\
&=\sum_{i,j=1}^n|\nabla
u(x_0)|^{p(x_0)-2}\Big(\delta_{ij}+(p(x_0)-2)\frac{{u}_{x_i}(x_0){u}_{x_j}(x_0)}{|\nabla
u(x_0)|^2}\Big){u}_{x_ix_j}(x_0)\\
&\qquad+\sum_{i=1}^n|\nabla
u(x_0)|^{p(x_0)-2}\Big({p}_{x_i}(x_0)\log|\nabla u(x_0)|\Big){u}_{x_i}(x_0) \\
&\qquad\qquad\ge \Delta_{p(x_0)}v (x_0).
\end{aligned}
\end{equation*}
That is, \eqref{u-supers} holds.

\bigskip

{\it Step II.} We now assume that $f$ and $p$ are as in the statement and we will show that $u$ is a viscosity supersolution to $\Delta_{p(x)}u=f$ in $\Omega$.

Again, let $v\in C^2(\Omega)$ such that  $v$ touches $u$ from below at $x_0\in \Omega$, with $\nabla v(x_0)\not=0$.
We will show that 
\begin{equation}\label{u-supersII}
\Delta_{p(x_0)}v(x_0)\leq f(x_0).
\end{equation}
Assume that $\Delta_{p(x_0)}v(x_0)> f(x_0)$. Then, there exist $r>0$ and $\sigma>0$ small such that
\begin{equation}\label{vgreatf}
\begin{aligned}
|\nabla v(x)|&> \sigma \quad \text{ in } B_r(x_0),\\
\Delta_{p(x)}v(x)&> f(x)+\sigma \quad \text{ in } B_r(x_0).
\end{aligned}
\end{equation}

We now take $p_k\in C^{1,\beta}(\overline{B_r(x_0)})$, for some $0<\beta<1$, with $\frac{1}{2}({1+p_{\min}})\le p_k(x)\le p_{\max}$, $p_k\le p$ in $B_r(x_0)$ and 
$\|\nabla p_k\|_{L^{\infty}}\leq 2L$, and
$f_k\in W^{1,{\infty}}(B_r(x_0))$, $||f_k||_{L^{\infty}}\le 2||f||_{L^{\infty}}$, such that
\begin{equation}\label{fkpkunif}
\begin{aligned}
&f_k\to f\quad \text{uniformly on } \overline{B_r(x_0)},\\
p_k\to p \  &\text{ and } \ \nabla p_k\to \nabla p \quad \text{uniformly on } \overline{B_r(x_0)}.
\end{aligned}
\end{equation}
Let $u_k\in W^{1,p_k(\cdot)}(B_r(x_0))$ be the (weak) solutions to 
\begin{equation*}
\begin{aligned}
\Delta_{p_k(x)}u_k&=f_k \text{ in }B_r(x_0),\\
u_k&=u \text{ on }\partial B_r(x_0).
\end{aligned}
\end{equation*}

Using Theorem 4.1 in \cite{FanZ}  and Theorem 1.2 in \cite{Fan}, we get that $u_k\in C^{1,\alpha}$ in  $\overline{B_r(x_0)}$, for some $0<\alpha<1$, $||u_k||_{C^{1,\alpha}(\overline{B_r(x_0)})}\le C$ and
\begin{equation}\label{conv-uk}
 u_k\to u\quad \text{uniformly on } \overline{B_r(x_0)}.
\end{equation}

Moreover, from the results in {\it Step I} we know that, for every $k$, $u_k$ is a viscosity supersolution to $\Delta_{p_k(x)}u_k=f_k$ in ${B_r(x_0)}$.

We fix $\ep>0$ and define
$$\widetilde v(x)= v(x)-{\ep}|x-x_0|^2.$$
Since there holds \eqref{vgreatf}, we can choose $\ep$ small enough so that 
\begin{equation}\label{tildevgreatf}
\begin{aligned}
|\nabla \widetilde v(x)|&> \frac{\sigma}{2} \quad \text{ in } B_r(x_0),\\
\Delta_{p(x)}\widetilde v(x)&> f(x)+\frac{\sigma}{2} \quad \text{ in } B_r(x_0).
\end{aligned}
\end{equation}

 Now, from
\eqref{fkpkunif} and \eqref{tildevgreatf}, we get
\begin{equation}\label{tildevgreatfk}
\Delta_{p_k(x)}\widetilde v(x)> f_k(x)+\frac{\sigma}{4} \quad \text{ in } B_r(x_0), \quad \text{ if }k\ge k_0.
\end{equation}
We now take $0<\delta<\frac{\ep}{4}r^2$. Recalling \eqref{conv-uk}, we can choose $k\ge k_0$ such that 
\begin{equation*}
|u_k-u|<\delta \quad \text{ in } \overline{B_r(x_0)},
\end{equation*}
so that we have
\begin{equation*}\label{uktildev}
\begin{aligned}
u_k+\delta & > \widetilde v \quad \text{ in } \overline{B_r(x_0)},\\
u_k(x_0)-\delta &< \widetilde v(x_0).
\end{aligned}
\end{equation*}
We now take 
$$
\bar{t}=\inf\Big\{ t\in\R \ / \ u_k+t\ge \widetilde v \quad \text{ in }\overline{B_r(x_0)}\Big\}.
$$
Then, $|\bar{t}|\le \delta$ and
\begin{equation}\label{ukbart}
\begin{aligned}
u_k & \ge \widetilde v -\bar{t}\quad \text{ in } \overline{B_r(x_0)},\\
u_k(\bar{x})&=\widetilde v(\bar{x})-\bar{t}, \quad\text{ for some } \bar{x}\in \overline{B_r(x_0)}.
\end{aligned}
\end{equation}
Suppose $\bar{x}\in  \partial {B_r(x_0)}$. Then,
$$u_k(\bar{x})=\widetilde v(\bar{x})-\bar{t}=v(\bar{x})-\ep r^2-\bar{t}\le u(\bar{x}) -\ep r^2+\delta\le u_k(\bar{x})+2\delta -\ep r^2,$$
a contradiction since we have chosen $\delta<\frac{\ep}{4}r^2$. 

Then $\bar{x}\in {B_r(x_0)}$ and \eqref{ukbart} says that  $\widetilde v-\bar{t}$ touches $u_k$ from below at $\bar{x}$. 
Since $\nabla\widetilde v(\bar{x})\neq0$, we get
\begin{equation*}
\Delta_{p_k(\bar{x})}\widetilde v(\bar{x})\le f_k(\bar{x}).
\end{equation*}
This contradicts \eqref{tildevgreatfk} and we conclude that   \eqref{u-supersII} holds. So $u$ is a viscosity supersolution to $\Delta_{p(x)}u=f$
in $\Omega$. 

\medskip

The proof that $u$ is a viscosity subsolution to $\Delta_{p(x)}u=f$ in $\Omega$ follows similarly.
\end{proof}

\begin{rem} \label{equiv-not}
As already mentioned in the Introduction,
the equivalence between weak and viscosity solutions to the $p(x)$-Laplacian with right hand side $f\equiv 0$ was proved in \cite{JLP}. On the other hand, this equivalence, in case $p(x)\equiv p$ and $f\not\equiv 0$ was dealt with in \cite{JJ} and \cite{MO}. See also \cite{JLM} for the case $p(x)\equiv p$ and $f\equiv 0$.
\end{rem}

\smallskip

We also obtain the following result that will be used in the proof of Lemma \ref{pre-harnack}

\begin{prop} \label{weak-is-strong} Let $p$ and $f$ be as in Definition \ref{defnweak}. Let $B_{2r}(x_0)\subset\subset\Omega$.

 Let $u\in W^{1,p(\cdot)}(\Omega)\cap L^{\infty}(\Omega)$ be a weak solution to $\Delta_{p(x)}u=f$ in $\Omega$ such that
$$c_1\le |\nabla u(x)| \le C_1 \quad \text{ in } B_{2r}(x_0), \quad c_1, C_1 \text{ positive constants}.$$

Then, $u\in W^{2,n}(B_r(x_0))$ and it is a strong solution to the linear uniformly elliptic equation
\[
\sum_{i,j=1}^n a_{ij}(x){u}_{x_ix_j}+\sum_{i=1}^n b_i(x){u}_{x_i}=f \quad \text{ in } B_r(x_0)
\]
where
$$a_{ij}(x)=|\nabla
u|^{p(x)-2}\Big(\delta_{ij}+(p(x)-2)\frac{{u}_{x_i}{u}_{x_j}}{|\nabla
u|^2}\Big),$$
and
$$b_i(x)=|\nabla
u|^{p(x)-2}\Big({p}_{x_i}(x)\log|\nabla u|\Big),$$ 
with
\begin{equation*}
\beta_1|\xi|^2\le \sum_{i,j=1}^n a_{ij}(x)\xi_i\xi_j\le
\beta_2|\xi|^2, \quad\forall \xi\in\R^N, \ \forall x\in B_r(x_0),
\end{equation*}
for $\beta_1, \beta_2$ positive constants, depending only on $c_1, C_1, p_{\min},  p_{\max}$. 
\end{prop}
\begin{proof}
We  take $f_k\in W^{1,{\infty}}(B_{2r}(x_0))$, $||f_k||_{L^{\infty}}\le 2||f||_{L^{\infty}}$, such that
\begin{equation*}
f_k\to f\quad \text{in } L^1({B_{2r}(x_0)}).
\end{equation*}
Let $u_k\in W^{1,p(\cdot)}(B_{2r}(x_0))$ be the (weak) solutions to 
\begin{equation*}
\begin{aligned}
\Delta_{p(x)}u_k&=f_k \text{ in }B_{2r}(x_0),\\
u_k&=u \text{ on }\partial B_{2r}(x_0).
\end{aligned}
\end{equation*}

Using Theorem 4.1 in \cite{FanZ}  and Theorem 1.2 in \cite{Fan}, we get that $u_k\in C^{1,\alpha}$ in  $\overline{B_{2r}(x_0)}$, for some $0<\alpha<1$, $||u_k||_{C^{1,\alpha}(\overline{B_{2r}(x_0)})}\le C$ and
\begin{equation*}
 u_k\to u, \quad  \nabla u_k\to \nabla u \quad\text{uniformly on } \overline{B_{2r}(x_0)}.
\end{equation*}
Then, for $k$ large,
$$\frac{c_1}{2}\le |\nabla u_k(x)| \le 2C_1 \quad \text{ in } B_{2r}(x_0).$$
Now, arguing as in Theorem 3.2 in \cite{CL}, we deduce that, for $k$ large, $u_k\in W^{2,2}_{\rm loc}(B_{2r}(x_0))$ and it is a solution to the linear uniformly elliptic equation
\[
\sum_{i,j=1}^n a^k_{ij}(x){(u_k)}_{x_ix_j}+\sum_{i=1}^n b^k_i(x){(u_k)}_{x_i}=f_k \quad \text{ in } B_{2r}(x_0)
\]
where
$$a^k_{ij}(x)=|\nabla
u_k|^{p(x)-2}\Big(\delta_{ij}+(p(x)-2)\frac{{(u_k)}_{x_i}{(u_k)}_{x_j}}{|\nabla
u_k|^2}\Big),$$
and
$$b^k_i(x)=|\nabla
u_k|^{p(x)-2}\Big({p}_{x_i}(x)\log|\nabla u_k|\Big),$$ 
with
\begin{equation*}
\beta_1|\xi|^2\le \sum_{i,j=1}^n a^k_{ij}(x)\xi_i\xi_j\le
\beta_2|\xi|^2, \quad\forall \xi\in\R^N, \ \forall x\in B_{2r}(x_0),
\end{equation*}
for $\beta_1, \beta_2$ positive constants, depending only on $c_1, C_1, p_{\min},  p_{\max}$. Moreover, 
$a^k_{ij}\in C^{\alpha}(\overline{B_{2r}(x_0)})$ and 
$||b^k_i||_{L^{\infty}(B_{2r}(x_0))}\le \bar{C}$.

It follows (see, for instance, Lemma 9.16 and Theorem 9.11  in \cite{GT}) that 
\begin{equation*}
u_k\in W^{2,n}_{\rm loc}(B_{2r}(x_0))\cap L^{\infty}(B_{2r}(x_0)) \quad \text{and}\quad ||u_k||_{W^{2,n}(B_{r}(x_0))}\le \tilde{C},
\end{equation*}
for some positive constant $\tilde{C}$. Then, passing to the limit $k\to\infty$, we get the desired result.
\end{proof}

\section{Auxiliary results}\label{section4}

In this section we prove some  results that will be of use in our main theorem. Namely, a Harnack inequality for an auxiliary problem of $p(x)$-Laplacian type and an existence result of barrier functions for the $p(x)$-Laplacian operator.

\smallskip

In the next result we assume for simplicity that $||f||_{L^{\infty}(\Omega)}\le 1$, but a similar result holds for any $f\in L^{\infty}(\Omega)$. We have

\begin{lem}\label{harnack-with-e} Assume that $1<p_{\min}\le p(x)\le p_{\max}<\infty$
with  $p(x)$ Lipschitz continuous in $\Omega$ and $\|\nabla p\|_{L^{\infty}}\leq L$, for some $L>0$.
Let $x_0\in\Omega$
and $0<R\le 1$ such that $\overline{B_{4R}(x_0)}\subset\Omega$.  Let $v\in W^{1,p(\cdot)}(\Omega)\cap L^{\infty}(\Omega)$
be a nonnegative solution to
\begin{equation}\label{eq-with-e}
\mbox{\rm div} (|\nabla v+e|^{p(x)-2}(\nabla v+e)) =f \quad\mbox{
in }\Omega,
\end{equation}
where $f\in L^{\infty}(\Omega)$ with $||f||_{L^{\infty}(\Omega)}\le 1$ and $e\in \mathbb{R}^n$ with $|e|=1$.  Then, there exists  $C$ such that
\begin{equation}\label{quasi_Harnack}\sup_{{B_R}(x_0)}v\leq C\Big[\inf_{B_R(x_0)}v+ R\Big({||f||_{L^{\infty}(B_{4R}(x_0))}}^{\frac{1}{p_{\max}-1}}+C\Big)\Big].\end{equation}
The constant $C$ depends only on $n$, $p_{\min}$,
$p_{\max}$, $||v||_{L^{\infty}(B_{4R}(x_0))}$  and $L$.
\end{lem}
\begin{proof}
We define $A:\Omega\times\mathbb{R}^n\to\mathbb{R}^n$
\begin{equation*}
A(x,\xi)=|\xi+e|^{p(x)-2}(\xi+e).
\end{equation*}
Then equation \eqref{eq-with-e} takes the form
\begin{equation*}
\mbox{\rm div}\, A(x,\nabla v) =f(x) \quad\mbox{
in }\Omega.
\end{equation*}

We first observe that, for  every $\xi\in \mathbb{R}^n$, 
\begin{equation*}
|A(x,\xi)|=|\xi+e|^{p(x)-1}\le C_1 |\xi|^{p(x)-1} +C_1,
\end{equation*}
where $C_1$ depends only on $p_{\max}$.
On the other hand, for  every $\xi\in \mathbb{R}^n$, 
\begin{equation}\label{lower-bound-A1}
\begin{aligned}
\langle A(x,\xi),\xi\rangle=&\:|\xi+e|^{p(x)-2}\langle\xi+e,\xi\rangle\\
=&\:|\xi+e|^{p(x)}-|\xi+e|^{p(x)-2}\langle\xi+e,e\rangle\\
\ge&\: |\xi+e|^{p(x)}-|\xi+e|^{p(x)-1}.
\end{aligned}
\end{equation}
Now, if $|\xi+e|\le 2$, we get from \eqref{lower-bound-A1}
\begin{equation}\label{lower-bound-A2}
\begin{aligned}
\langle A(x,\xi),\xi\rangle\ge&\:|\xi+e|^{p(x)}-2^{p(x)-1}\\
\ge&\: C_2|\xi|^{p(x)}-C_3,
\end{aligned}
\end{equation}
where $C_2$ and $C_3$ depend only on $p_{\max}$. If $|\xi+e|> 2$, we obtain from \eqref{lower-bound-A1}
\begin{equation}\label{lower-bound-A3}
\begin{aligned}
\langle A(x,\xi),\xi\rangle\ge&\:|\xi+e|^{p(x)}-|\xi+e|^{p(x)-1}\\
=&\:|\xi+e|^{p(x)}(1-|\xi+e|^{-1})\\
\ge&\:\frac12 |\xi+e|^{p(x)}
\ge C_4|\xi|^{p(x)}-\frac12,
\end{aligned}
\end{equation}
where $C_4$ depends only on $p_{\max}$. Then, from \eqref{lower-bound-A2} and \eqref{lower-bound-A3} we deduce
\begin{equation*}
\langle A(x,\xi),\xi\rangle\ge C_5|\xi|^{p(x)}-C_6,
\end{equation*}
where $C_5$ and $C_6$ depend only on $p_{\max}$. Now the result follows from Theorem 1.1 in \cite{Wo}.
\end{proof}

\smallskip

We now continue with a technical result concerning the existence of barrier functions  for the $p(x)$-Laplacian operator.

\begin{lem}\label{barry} Let $x_0\in B_1$ and $0<\bar{r}_1<\bar{r}_2\le 1$. 
Assume that $1<p_{\min}\le p(x)\le p_{\max}<\infty$
 and $\|\nabla p\|_{L^{\infty}}\leq \ep^{1+\theta}$, for some $0<\theta\le 1$. Let $c_0, c_1, c_2$ be positive constants and let and
$c_3\in \R$.

There exist positive constants $\gamma\ge 1$, $\bar{c}$, $\ep_0$ and $\ep_1$ such that the functions
\begin{equation*}
w(x)=c_1|x-x_0|^{-\gamma}-c_2, 
\end{equation*}
\begin{equation*}
v(x)=q(x)+\frac{c_0}{2}\varepsilon (w(x)-1),\quad q(x)=x_n+c_3
\end{equation*}
satisfy, for $\bar{r}_1\le |x-x_0|\le \bar{r}_2$,
\begin{equation}\label{eq-w}
\Delta_{p(x)}w\geq \bar{c}, \quad \text{for }\, 0<\ep\le \ep_0,
\end{equation}
\begin{equation}\label{eq-v}
\frac{1}{2}\le |\nabla v|\le 2, \qquad \Delta_{p(x)}v > \varepsilon^2, \quad \text{for }\, 0<\ep\le \ep_1.
\end{equation}
Here $\gamma=\gamma(n,p_{\min}, p_{\max})$, $\bar{c}=\bar{c}(p_{\min}, p_{\max},  c_1)$, $\ep_0=\ep_0(n,p_{\min}, p_{\max}, \bar{r}_1,  c_1)$, $\ep_1=\ep_1(n,p_{\min}, p_{\max}, \bar{r}_1,  c_0, c_1, \theta)$.
\end{lem}
\begin{proof} Without loss of generality we can assume that $x_0=0$. We will divide the proof into five steps.

 \vspace{2mm}

{\it Step 1.} For simplicity, we assume first that $c_1=1$. Let us fix $p\in\R$, $1<p_{\min}\le p\le p_{\max}<\infty$ and $\gamma>0$. Let us consider $x\in{{\R}^n}\setminus\{0\}$.

Then, $w(x)=|x|^{-\gamma}-c_2$ and $\nabla w=-\gamma |x|^{-\gamma-2}x$, so that 
$$
\frac{\nabla w}{|\nabla w|}=-\frac{x}{|x|}.
$$
Moreover
\begin{equation}\begin{split}
D^2 w&=\gamma (\gamma+2)|x|^{-\gamma-2}\frac{x}{|x|}\otimes \frac{x}{|x|}-\gamma |x|^{-\gamma-2}I\\
&=\gamma |x|^{-\gamma-2}\left((\gamma+2)\frac{x}{|x|}\otimes \frac{x}{|x|}-I\right).
\end{split}
\end{equation}
As a consequence
\begin{equation}\begin{split}
\mbox{Tr}(D^2 w)=\gamma |x|^{-\gamma-2}\left((\gamma+2)-n\right).\end{split}
\end{equation}
Thus
\begin{equation}\label{421}\begin{split}
&\Delta_p w=|\nabla w|^{p-2}\left(\Delta w+(p-2)\langle D^2w\frac{\nabla w}{|\nabla w|}, \frac{\nabla w}{|\nabla w|}\rangle\right)\\
&=\gamma^{p-1} |x|^{-(\gamma+1)(p-2)} |x|^{-\gamma-2}\left((\gamma+2)-n+(p-2)\langle[(\gamma+2)\frac{x}{|x|}\otimes \frac{x}{|x|}-I]\frac{x}{|x|},\frac{x}{|x|}\rangle\right)\\
&=\gamma^{p-1} |x|^{-\gamma(p-1)-p}\left(\gamma+2-n+(p-2)(\gamma+1)\right)\\
&=\gamma^{p-1} |x|^{-\gamma(p-1)-p}\left(\gamma(p-1)+p-n\right)\\
&\ge \gamma^{p-1} |x|^{-\gamma(p-1)-p}\left(\gamma(p_{\min}-1)+p_{\min}-n\right)\ge  \gamma^{p-1} |x|^{-\gamma(p-1)-p},
\end{split}
\end{equation}
if $\gamma>0$ is such that
\begin{equation}\label{cond-1}
\gamma(p_{\min}-1)+p_{\min}-n \ge 1.
\end{equation}

On the other hand,
$$
D^2v=\frac{c_0}{2}\varepsilon D^2 w.
$$
Then, for $x$ such that $\nabla v(x)\neq 0$,
\begin{equation}\label{pilu}
\begin{split}
&\Delta_p v=|\nabla v|^{p-2}\left(\Delta v+(p-2)\langle D^2v\frac{\nabla v}{|\nabla v|}, \frac{\nabla v}{|\nabla v|}\rangle\right)\\
&=\frac{c_0}{2}\varepsilon|\nabla v|^{p-2}\left(\Delta w+(p-2)\langle D^2w\frac{\nabla v}{|\nabla v|}, \frac{\nabla v}{|\nabla v|}\rangle\right)\\
&=\frac{c_0}{2}\gamma\varepsilon|\nabla v|^{p-2} |x|^{-\gamma-2}\left\{(\gamma+2)-n+(p-2)\langle[(\gamma+2)\frac{x}{|x|}\otimes \frac{x}{|x|}-I]\frac{\nabla v}{|\nabla v|},\frac{\nabla v}{|\nabla v|}\rangle\right\}\\
&=\frac{c_0}{2}\gamma\varepsilon|\nabla v|^{p-2} |x|^{-\gamma-2}\left\{(\gamma+2)-n+(p-2)\left[(\gamma+2)\langle\frac{x}{|x|},\frac{\nabla v}{|\nabla v|}\rangle^2-1\right]\right\}\\
&=\frac{c_0}{2}\gamma\varepsilon|\nabla v|^{p-2} |x|^{-\gamma-2}\left\{(\gamma+2)\left[1+(p-2)\langle\frac{x}{|x|},\frac{\nabla v}{|\nabla v|}\rangle^2\right]-n-p+2\right\}.\\
\end{split}
\end{equation}
We also observe that 
\begin{equation}
0\le \langle\frac{x}{|x|},\frac{\nabla v}{|\nabla v|}\rangle^2\le 1.
\end{equation}
Hence, in case $p_{\min}\le p\leq 2,$   it follows from \eqref{pilu}
\begin{equation}\label{pilu1}
\begin{split}
\Delta_p v&\geq \frac{c_0}{2}\gamma\varepsilon|\nabla v|^{p-2} |x|^{-\gamma-2}\left\{(\gamma+2)(1+p-2)-n-p+2\right\}\\
&=\frac{c_0}{2}\gamma\varepsilon|\nabla v|^{p-2} |x|^{-\gamma-2}\left\{(\gamma+2)(p-1)-n-p+2\right\}\\
&\ge\frac{c_0}{2}\gamma\varepsilon|\nabla v|^{p-2} |x|^{-\gamma-2}\left\{(\gamma+2)(p_{\min}-1)-n\right\}
\ge \frac{c_0}{2}\gamma\varepsilon|\nabla v|^{p-2} |x|^{-\gamma-2},
\end{split}
\end{equation}
if $\gamma>0$ is such that
\begin{equation}\label{cond-2}
(\gamma+2)(p_{\min}-1)-n \ge 1.
\end{equation}

Moreover, in case $2<p\le p_{\max},$ it follows from \eqref{pilu} 
\begin{equation}\label{pilu2}
\begin{split}
\Delta_pv&\geq \frac{c_0}{2}\gamma\varepsilon|\nabla v|^{p-2} |x|^{-\gamma-2}\left(\gamma+2-n-p+2\right)\\
&\ge\frac{c_0}{2}\gamma\varepsilon|\nabla v|^{p-2} |x|^{-\gamma-2}\left(\gamma+4-n-p_{\max}\right)
\ge\frac{c_0}{2}\gamma\varepsilon|\nabla v|^{p-2} |x|^{-\gamma-2},
\end{split}
\end{equation}
if $\gamma>0$ is such that
\begin{equation}\label{cond-3}
\gamma+4-n-p_{\max} \ge 1.
\end{equation}

We now fix 
 \begin{equation}\label{fix-gamma}
\gamma=\gamma(n,p_{\min},p_{\max})=\max\left\{1, \ \ \frac{1+n-p_{\min}}{p_{\min}-1}, \ \ \frac{1+n}{p_{\min}-1}-2, \ \ n+p_{\max}-3\right\}.
\end{equation}

Then, $\gamma=\gamma(n,p_{\min},p_{\max})\ge 1$ and $\gamma$ satisfies  \eqref{cond-1}, \eqref{cond-2} and \eqref{cond-3}. Hence we obtain 
 from \eqref{421}, \eqref{pilu1} and \eqref{pilu2} that for every $p\in [p_{\min},p_{\max}]$ and $x\in{{\R}^n}\setminus\{0\}$
\begin{equation}\label{421b}
\Delta_p w \ge |x|^{-\gamma(p-1)-p},
\end{equation}
\begin{equation}\label{xxplus}
\Delta_{p}v\geq \frac{c_0}{2}\varepsilon|\nabla v|^{p-2} |x|^{-\gamma-2},\quad \text{if } \nabla v(x)\neq 0.
\end{equation}

\vspace{2mm}

{\it Step 2.} We now assume that  $c_1>0$ is arbitrary. We fix $\gamma=\gamma(n,p_{\min},p_{\max})\ge 1$ as above,  given by 
\eqref{fix-gamma}. It is not hard to see that similar  computations as those in Step 1, but with $c_1>0$ arbitrary, imply that 
for every $p\in [p_{\min},p_{\max}]$ and $x\in{{\R}^n}\setminus\{0\}$
\begin{equation}\label{421b-with-c1}
\Delta_p w \ge {c_1}^{p-1}|x|^{-\gamma(p-1)-p},
\end{equation}
\begin{equation}\label{xxplus-with-c1}
\Delta_{p}v\geq \frac{c_0}{2}c_1\varepsilon|\nabla v|^{p-2} |x|^{-\gamma-2},\quad \text{if } \nabla v(x)\neq 0.
\end{equation}

\vspace{2mm}

{\it Step 3.} We now observe that there holds 
$$
\nabla v=e_n+\frac{c_0}{2}\varepsilon \nabla w.
$$
Then, for $\bar{r}_1\le |x|\le \bar{r}_2$,
\begin{equation*}
\begin{aligned}
\Big| |\nabla v|-1\Big|&= \Big| |\nabla v|-|e_n|\Big|\le \Big| \nabla v-e_n\Big|=\Big|\frac{c_0}{2}\varepsilon \nabla w\Big|\\
&=
\frac{c_0}{2}c_1\varepsilon \gamma |x|^{-\gamma-1}\le \frac{c_0}{2}c_1\varepsilon \gamma \bar{r}_1^{-\gamma-1}\le \frac{1}{2},
\end{aligned}
\end{equation*}
if we let $\ep\le \bar{\ep}_1=\bar{\ep}_1(n,p_{\min},p_{\max}, \bar{r}_1, c_0, c_1)$ and therefore,
\begin{equation}\label{eq-v-1}
\frac{1}{2}\le |\nabla v|\le 2, \quad \text{for }\, \ep\le \bar{\ep}_1.
\end{equation}
 So the first assertion in \eqref{eq-v} follows.

\vspace{2mm}

{\it Step 4.} We now consider $p(x)$ a Lipschitz continuous function such that $1<p_{\min}\le p(x)\le p_{\max}<\infty$.

We first observe that, for any $R\ge 1$, 
\begin{equation*}
\begin{aligned}
& t^{p(x)-1}\big|\log t\big|\le t^{p_{\min}-1}\big|\log t\big|\le C_1(p_{\min}), \qquad \text{if }\, 0<t<1,\\
& t^{p(x)-1}\big|\log t\big|\le t^{p_{\max}-1}\big|\log t\big|\le R^{p_{\max}-1} \log R \qquad \text{if }\, 1\le t\le R,
\end{aligned}
\end{equation*}
so that
\begin{equation}\label{bound-for-log}
t^{p(x)-1}\big|\log t\big|\le C_2(p_{\min},p_{\max}, R), \qquad \text{if }\, 0<t\le R.
\end{equation}
It then follows from \eqref{eq-v-1} and \eqref{bound-for-log} that, for $\bar{r}_1\le |x|\le \bar{r}_2$,
\begin{equation}\label{bound-nabla v}
|\nabla v|^{p(x)-1}\big|\log |\nabla v| \big|\le C_3(p_{\min},p_{\max}), \qquad \text{if }\,\ep\le \bar{\ep}_1.
\end{equation}
We also have, for $\bar{r}_1\le |x|\le \bar{r}_2$,
$$|\nabla w|=c_1\gamma |x|^{-\gamma-1}\le c_1\gamma \bar{r}_1^{-\gamma-1},
$$
so using again \eqref{bound-for-log}, we get, for $\bar{r}_1\le |x|\le \bar{r}_2$,
\begin{equation}\label{bound-nabla w}
|\nabla w|^{p(x)-1}\big|\log |\nabla w| \big|\le C_4(n,p_{\min},p_{\max}, \bar{r}_1, c_1).
\end{equation}

\vspace{2mm}

{\it Step 5.} We now assume that  $p(x)$  satisfies moreover that $\|\nabla p\|_{L^{\infty}}\leq \ep^{1+\theta}$, for some $0<\theta<1$. Then,
from \eqref{bound-nabla v} we obtain, for $\bar{r}_1\le |x|\le \bar{r}_2$,
\begin{equation}\label{bound-log-part}
\Big||\nabla v|^{p(x)-2}\langle \nabla p(x),\nabla v\rangle\log|\nabla v| \Big|\le |\nabla v|^{p(x)-1}\big|\log |\nabla v| \big|\|\nabla p\|_{L^{\infty}}\le \ep^{1+\theta} C_3,
\end{equation}
if $\ep\le \bar{\ep}_1$. 
Hence, from \eqref{xxplus-with-c1},\eqref{bound-log-part} and \eqref{eq-v-1}, for $\bar{r}_1\le |x|\le \bar{r}_2$,
\begin{equation}\begin{split}
\Delta_{p(x)}v&=|\nabla v|^{p(x)-2}(\Delta v+(p(x)-2)\langle D^2v\frac{\nabla v}{|\nabla v|},\frac{\nabla v}{|\nabla v|}\rangle+\langle \nabla p(x),\nabla v\rangle\log|\nabla v|\rangle)\\
&\geq \frac{c_0}{2}c_1\varepsilon|\nabla v|^{p(x)-2} |x|^{-\gamma-2}-\ep^{1+\theta} C_3 \\
&\ge \frac{c_0}{2}c_1\varepsilon C_5 |x|^{-\gamma-2}-\ep^{1+\theta} C_3\ge \frac{c_0}{2}c_1\varepsilon C_5-\ep^{1+\theta} C_3
=\varepsilon(\frac{c_0}{2}c_1C_5-\ep^{\theta} C_3),\\
\end{split}
\end{equation}
if $\ep\le \bar{\ep}_1$, where we have used that $\bar{r}_2\le 1$ and $C_5=C_5(p_{\min}, p_{\max})$, $C_5= \min\{ (\frac{1}{2})^{p_{\max}-2}, 2^{p_{\min}-2}\}$.
We conclude that, for $\bar{r}_1\le |x|\le \bar{r}_2$, 
\begin{equation*}
\Delta_{p(x)}v \ge \varepsilon(\frac{c_0}{2}c_1C_5-\ep^{\theta} C_3)\ge \varepsilon\frac{c_0}{4}c_1C_5>{\ep}^2,
\end{equation*}
if moreover $\ep\le\tilde{\ep}_1= \tilde{\ep}_1(p_{\min}, p_{\max}, c_0, c_1, \theta)$. That is, the second assertion in  \eqref{eq-v} follows.

Finally, from \eqref{bound-nabla w} we obtain, for $\bar{r}_1\le |x|\le \bar{r}_2$,
\begin{equation}\label{bound-log-part-w}
\Big||\nabla w|^{p(x)-2}\langle \nabla p(x),\nabla w\rangle\log|\nabla w| \Big|\le |\nabla w|^{p(x)-1}\big|\log |\nabla w| \big|\|\nabla p\|_{L^{\infty}}\le \ep^{1+\theta} C_4.
\end{equation}
Hence, from \eqref{421b-with-c1} and \eqref{bound-log-part-w}, for $\bar{r}_1\le |x|\le \bar{r}_2$,
\begin{equation}\begin{split}
\Delta_{p(x)}w&=|\nabla w|^{p(x)-2}(\Delta w+(p(x)-2)\langle D^2w\frac{\nabla w}{|\nabla w|},\frac{\nabla w}{|\nabla w|}\rangle+\langle \nabla p(x),\nabla w\rangle\log|\nabla w|\rangle)\\
&\geq c_1^{p(x)-1}|x|^{-\gamma(p(x)-1)-p(x)}-\ep^{1+\theta} C_4 \ge 2\bar{c}-\ep C_4,\\
\end{split}
\end{equation}
if $\ep\le 1$. Here we have used that $\bar{r}_2\le 1$ and we have denoted 
$\bar{c}=\bar{c}(p_{\min}, p_{\max}, c_1)=\frac{1}{2}\min\{ {c_1}^{p_{\min}-1},{c_1}^{p_{\max}-1}\}$. We conclude that, for $\bar{r}_1\le |x|\le \bar{r}_2$,
\begin{equation*}
\Delta_{p(x)}w \ge \bar{c},
\end{equation*}
if $\ep\le{\ep}_0= {\ep}_0(n,p_{\min}, p_{\max}, \bar{r}_1, c_1)$. This proves \eqref{eq-w} and finishes the proof.
\end{proof}

\section{Geometric regularity results}\label{section5}
 In this section we  prove a Harnack type inequality for a solution $u$ to problem \eqref{fb}, following the approach in \cite{D}.
We will argue  assuming that 
\begin{equation}\label{assumptions-harn}
||f||_{L^{\infty}(\Omega)}\le \ep^2,\quad ||g-1||_{L^{\infty}(\Omega)}\le \ep^2, \quad  ||\nabla p||_{L^{\infty}(\Omega)}\le \ep^{1+\theta},
\quad ||p-p_0||_{L^{\infty}(\Omega)}\le \ep,
\end{equation}
holds, for $0<\ep<1$, for some constant $0<\theta\le$1.

\smallskip

The proof of Harnack inequality is based on the following lemma.

\begin{lem}\label{pre-harnack}
Let $u$ be a  solution to  \eqref{fb}--\eqref{assumptions-harn} in $B_1$. There exists a universal constant $\bar{\varepsilon}$ such that if $0<\varepsilon \leq \bar{\varepsilon}$ and $u$ satisfies
\begin{equation}\label{flat}
q^+(x)\leq u(x)\leq (q(x)+\varepsilon)^+,\quad x\in B_1, \ \ q(x)=x_n+\sigma, \ \ |\sigma|<\frac{1}{20},
\end{equation} 
and in $x_0=\frac{1}{10}e_n,$ 
$$
u(x_0)\geq (q(x_0)+\frac{\varepsilon}{2})^+,
$$
then
\begin{equation}
 u \ge (q+c\varepsilon)^+ \quad \text{in } \overline{B}_{\frac{1}{2}},
\end{equation}
for some universal $0<c<1.$ Analogously, if 
\begin{equation}
u(x_0)\leq (q(x_0)+\frac{\varepsilon}{2})^+,
\end{equation}
then 
\begin{equation}
u\leq (q+(1-c)\varepsilon)^+\quad \text{in } \overline{B}_{\frac{1}{2}}.
\end{equation}
\end{lem}
\begin{proof}
The proof  follows the original one in \cite{D} adapted by the dichotomy discussed in \cite{LR}. We will prove the first statement. 

{}From \eqref{flat} we have that $u\ge q$ in $B_1$.

We also notice that $B_{1/20}(x_0)\subset B_1^+(u).$ Then, 
\begin{equation}\label{eq-p(x)}
\Delta_{p(x)}u=f\quad \text{ in } B_{1/20}(x_0).
\end{equation}
 Thus,  by Theorem 1.1 
in \cite{Fan}, $u\in C^{1,\alpha}$ in $\overline{B}_{1/40}(x_0),$  where $\alpha=\alpha (p_{\min},p_{\max},n)\in (0,1)$ and
$||u||_{C^{1,\alpha}(\overline{B}_{1/40}(x_0))}\leq C,$ with $C=C (p_{\min},p_{\max},n)\geq 1.$ Here we have used \eqref{assumptions-harn}  and also that \eqref{flat} implies that $||u||_{L^{\infty}(B_1)}\le 3$.

We will consider two cases:

\medskip

{\bf Case (i).} Suppose $|\nabla u(x_0)|<\frac{1}{4}$. We choose $r_1>0$, 
$r_1=r_1(p_{\min},p_{\max},n)\le 1/40$ such that $|\nabla u(x)|\le\frac{1}{2}$ in $B_{{r_1}}(x_0).$ In addition, there exists a constant  $0<r_2=r_2(r_1)=r_2(p_{\min},p_{\max},n)<r_1$ such that $(x-r_2e_n)\in B_{r_1}(x_0),$  for every $x\in B_{r_1/2}(x_0)$. We observe that $\tilde{v}=u-q$ satisfies
\begin{equation}
\mbox{\rm div} (|\nabla \tilde{v}+e_n|^{p(x)-2}(\nabla \tilde{v}+e_n)) =f \quad\mbox{
in } B_{\frac{1}{20}}(x_0).
\end{equation}
  We now apply Lemma \ref{harnack-with-e} to the function $\tilde{v}=u-q$  in $B_{4r_3}(x_0)$, 
 where $r_3=\min\{\frac{r_1}{4},\frac{r_2}{8}\}.$ 
 In particular we obtain from \eqref{quasi_Harnack} that
\begin{equation*}\begin{split}
u(x)-q(x)&\geq C^{-1}(u(x_0)-q(x_0))- r_3\geq\frac{\varepsilon}{2C}- r_3,
\end{split}
\end{equation*}
for  $x\in B_{r_3}(x_0)$. Here $C=C(n,p_{\min}, p_{\max})$ is a universal constant because 
$||f||_{L^{\infty}(B_1)}\le \ep^2,$ see \eqref{assumptions-harn}, and $||\tilde{v}||_{L^{\infty}(B_1)}\le 2$.

On the other hand, for all $x\in B_{r_3}(x_0)$ we obtain
\begin{equation*}\begin{split}
&\frac{\varepsilon}{2C}-  r_3\leq u(x)-q(x)=u((x-r_2e_n)+r_2e_n)-q((x-r_2e_n)+r_2e_n)\\
&=u((x-r_2e_n)+r_2e_n)-q(x-r_2e_n)-r_2\leq u(x-r_2e_n)-q(x-r_2e_n)+\frac{r_2}{2}-r_2.
\end{split}
\end{equation*}
As a consequence,   denoting $c_0=C^{-1}$ and $\bar{x}_0:=x_0-r_2e_n$, we get  for all $x\in B_{r_3}(\bar{x}_0)$
\begin{equation}\label{contracontra}\begin{split}
&\frac{c_0}{2}\varepsilon=\frac{\varepsilon}{2C}\leq\frac{\varepsilon}{2C}-  r_3+\frac{r_2}{2}=\frac{\varepsilon}{2C}- r_3-\frac{r_2}{2}+r_2\leq u(x)-q(x).
\end{split}
\end{equation}
Let us  define the function $w:\bar{D}\to\mathbb{R},$ $D:=B_{\frac{4}{5}}(\bar{x}_0)\setminus \bar{B}_{r_3}(\bar{x}_0)$ as
$$
w(x)=c\left(|x-\bar{x}_0|^{-\gamma}-(\frac{4}{5})^{-\gamma}\right),
$$
for $\gamma=\gamma(n,p_{\min}, p_{\max})\ge 1$ given in Lemma \ref{barry} (see \eqref{fix-gamma}). We choose $c=c(n,p_{\min}, p_{\max})>0$ in such a way that
$$
w=\left\{\begin{array}{l}
0,\quad \mbox{on}\quad \partial B_{\frac{4}{5}}(\bar{x}_0)\\
1,\quad \mbox{on}\quad \partial B_{r_3}(\bar{x}_0).
\end{array}
\right.
$$
As usual, we define for every $x\in \bar{B}_{\frac{4}{5}}(\bar{x}_0)$
$$
v(x)=q(x)+c_0\frac{\varepsilon}{2}(w(x)-1)
$$
and for $t\geq 0$ we set 
$$
v_t(x)=v(x)+t,\quad x\in\bar{B}_{\frac{4}{5}}(\bar{x}_0).
$$
We extend $w$ to $1$ in $B_{r_3}(\bar{x}_0),$ so that it results
$$
v_0(x)=v(x)\leq q(x)\leq u(x),\quad x\in \bar{B}_{\frac{4}{5}}(\bar{x}_0).
$$
Let
$$
\bar{t}=\sup\{t\geq 0:\:\:v_t\leq u\:\:\mbox{in}\:\:\bar{B}_{\frac{4}{5}}(\bar{x}_0)\}.
$$
{\bf Claim:} $\bar{t}\geq \frac{c_0\varepsilon}{2}.$

Assuming that the previous Claim holds, we obtain from the definition of $v$ that, in $B_{\frac{4}{5}}(\bar{x}_0),$ the  inequality
$$
u(x)\geq v(x)+\bar{t}\geq q(x)+\frac{c_0\varepsilon}{2}w(x)
$$
is satisfied.

On the other hand, $B_{\frac{1}{2}}\subset B_{\frac{3}{5}}(\bar{x}_0)$ and since
$$
w(x)\geq\left\{
\begin{array}{ll}
c\left((\frac{3}{5})^{-\gamma}-(\frac{4}{5})^{-\gamma}\right),& B_{\frac{3}{5}}(\bar{x}_0)\setminus B_{r_3}(\bar{x}_0),\\
1,& B_{r_3}(\bar{x}_0),
\end{array}
\right.
$$
we conclude that, in $B_{\frac{1}{2}}$, 
$$
u(x)-q(x)\geq c_1\varepsilon,
$$
with $0<c_1=c_1(n,p_{\min}, p_{\max})<1$ universal, as desired.

We now have to prove the Claim. We argue by contradiction assuming that 
$\bar{t}< \frac{c_0\varepsilon}{2}.$ Let $y_0\in \bar{B}_{\frac{4}{5}}(\bar{x}_0)$ be the contact point between $v_{\bar{t}}$ and $u$, where
$$
v_{\bar{t}}(y_0)=u(y_0).
$$ 

We will prove that $y_0\in \overline{B}_{r_3}(\bar{x}_0).$ In fact, recalling that $w$ vanishes on $\partial B_{\frac{4}{5}}(\bar{x}_0)$ and from the definition of $v_{\bar{t}}$, we obtain
$$
v_{\bar{t}}=q-\frac{c_0}{2}\varepsilon+\bar{t}<u\quad \mbox{on}\:\: \partial B_{\frac{4}{5}}(\bar{x}_0),
$$
because $u\geq q$ and $\bar{t}<\frac{c_0\varepsilon}{2}.$

We can apply Lemma \ref{barry} to $v.$ Hence, there exists $\varepsilon_1=\ep_1(n,p_{\min}, p_{\max}, \theta)$ a  universal constant such that 
$$
\frac{1}{2}\leq |\nabla v_{\bar{t}}|=|\nabla v|\leq 2,
$$
$$
\Delta_{p(x)}v_{\bar{t}}=\Delta_{p(x)}v> {\varepsilon}^2\ge f,
$$
 for every $0<\varepsilon\leq\varepsilon_1$ and for every $x\in D= B_{\frac{4}{5}}(\bar{x}_0)\setminus \overline{B}_{r_3}(\bar{x}_0).$ 

On the other hand, from the definition of  $v_{\bar{t}}$, we have
\begin{equation}\label{cfrf}
|\nabla v_{\bar{t}}|\geq |(v_{\bar{t}})_n|=|1+\frac{c_0}{2}\varepsilon w_n|,
\end{equation}
where $(v_{\bar{t}})_n$ and $w_n$ denote the partial derivatives with respect to $x_n$ of $v_{\bar{t}}$ and $w$. 

Let us show that $w_n>\hat{c}$ in $\{v_{\bar{t}}\leq 0\}\cap D$, for $\hat{c}>0$ universal. 

In fact, whenever $0<\varepsilon\le \ep_2$, for $\ep_2$ universal, we have 
 $$
\{v_{\bar{t}}\leq 0\}\cap D\subset\{q\leq \frac{c_0\varepsilon}{2}\}= \{x_n\leq \frac{c_0\varepsilon}{2}-\sigma\}\subset\{x_n\leq \frac{5}{80}\}.
$$
On the other hand,
$$
\nabla w=-\gamma c|x-\bar{x}_0|^{-\gamma-2}(x-\bar{x}_0)=-\gamma c|x-\bar{x}_0|^{-\gamma-1}\frac{x-\bar{x}_0}{|x-\bar{x}_0|}.
$$
Moreover, denoting $\nu_x=\frac{x-\bar{x}_0}{|x-\bar{x}_0|},$ we observe that, in $\{v_{\bar{t}}\leq 0\}\cap D$, we have
$-\langle \nu_x,e_n\rangle>0$ since 
$$x_n-(\bar{x}_0)_n=x_n-\frac{1}{10}+r_2\leq -\frac{1}{80}\quad \text{ in } \{x_n\leq \frac{5}{80}\}.$$ 
In particular, there holds in $\{v_{\bar{t}}\leq 0\}\cap D$
$$w_n=\langle\nabla w,e_n\rangle=-\gamma\langle \nu_x,e_n\rangle c|x-\bar{x}_0|^{-\gamma-1}\geq c\gamma\frac{1}{80}\frac{5}{4}(\frac{4}{5})^{-1-\gamma}=\hat{c}>0.$$

Thus, from \eqref{cfrf} we deduce that 
$$
|\nabla v_{\bar{t}}|\geq 1+\frac{c_0}{2}\varepsilon w_n\geq 1+\frac{c_0}{2}\hat{c}\varepsilon
$$
in $\{v_{\bar{t}}\leq 0\}\cap D$, which implies, for $\varepsilon$ sufficiently small,
$$
|\nabla v_{\bar{t}}|> 1+\varepsilon^2\geq g,
$$
on ${F}(v_{\bar{t}})\cap D.$ Then $v_{\bar{t}}$ is a strict subsolution to \eqref{fb} in $D$ touching $u$ at $y_0.$ Hence $y_0\in \overline{B}_{r_3}(\bar{x}_0)$ and this generates  a contradiction with \eqref{contracontra},  because 
$$
u(y_0)=v_{\bar{t}}(y_0)=v(y_0)+\bar{t} = q(y_0)+\bar{t}<q(y_0)+c_0\varepsilon.
$$

\medskip

{\bf Case (ii).} Now suppose $|\nabla u(x_0)|\geq\frac{1}{4}.$ By exploiting  the $C^{1,\alpha}$ regularity of $u$ in $\overline{B}_{\frac{1}{40}}(x_0)$, we know that $u$ is Lipschitz continuous in $\overline{B}_{\frac{1}{40}}(x_0)$, as well as there exist a  constant $0<r_0=r_0(n, p_{\min},p_{\max})$, with $8r_0\le \frac{1}{40}$, and $C=C(n,p_{\min},p_{\max})>1$ such that 
$$
\frac{1}{8}\leq |\nabla u|\leq C \quad \text{ in } B_{8r_0}(x_0).
$$ 
In addition, since \eqref{eq-p(x)} holds, it follows by Proposition \ref{weak-is-strong}, that
$u\in W^{2,n}(B_{4r_0}(x_0))$ and it is a solution to the linear uniformly elliptic 
equation
$$
\mathcal{L}h=f \quad \text{ in } B_{4r_0}(x_0),
$$
where
\begin{equation*}\begin{split}
&\mathcal{L}h=\mbox{Tr}(A(x)D^2h(x))+\langle b,\nabla h(x)\rangle,
\end{split}
\end{equation*}
$$
A(x):=|\nabla u|^{p(x)-2}\left(I+(p(x)-2)\frac{\nabla u(x)}{|\nabla u(x)|}\otimes \frac{\nabla u(x)}{|\nabla u(x)|}\right),
$$
and
$$
b(x):=|\nabla u|^{p(x)-2}\log|\nabla u(x)|\nabla p(x).
$$
Hence
$A\in C^{0,\alpha}(\overline{B}_{4r_0}(x_0))$, $b\in C(\overline{B}_{4r_0}(x_0))$   and  $\mathcal{L}$ has universal ellipticity constants (depending only on $n, p_{\min},p_{\max}$). Moreover, $||b||_{L^{\infty}(B_{4r_0}(x_0))}\leq C\varepsilon^{1+\theta}$, $C$ universal, because $ ||\nabla p||_{L^{\infty}(B_1)}\le \ep^{1+\theta}$ (see \eqref{assumptions-harn}).

In this way, we conclude that $u-q$ satisfies
\begin{equation*}\begin{split}
&\mbox{Tr}(A(x)D^2h(x))+\langle b,\nabla h(x)\rangle =f-\langle b,e_n\rangle \quad \text{ in } B_{4r_0}(x_0).
\end{split}
\end{equation*}
 Then,  applying  Harnack's inequality (see, for instance, \cite{GT}, Chap. 9) and recalling again \eqref{assumptions-harn}, we obtain
\begin{equation}
\begin{aligned}
u(x)-q(x)\geq  C_1(u(x_0)-&q(x_0))-C_2(||f||_{L^\infty(B_{4r_0}(x_0))}+||b||_{L^\infty(B_{4r_0}(x_0))})\\
  &\geq C_1\frac{\varepsilon}{2}-C_2(\ep^2+C\varepsilon^{1+\theta})\geq \frac{c_0}{2}\varepsilon,
\end{aligned}
\end{equation}
for every $x\in B_{r_0}(x_0)$, for $0<\varepsilon\le \ep_3$. Here $\ep_3$, $C_1$, $C_2$ and $c_0$ are positive universal constants. At this point, we can repeat the same argument of Case (i) around the point $x_0$, considering the annulus 
$B_{\frac{4}{5}}({x}_0)\setminus \bar{B}_{r_0}({x}_0)$. This completes the proof.
\end{proof}

The next  result is the main tool in Theorem \ref{flatmain1}.

\begin{thm}[Harnack inequality]
\label{HI}
There exists a universal constant $\bar
\ep$,  such that if $u$ solves \eqref{fb}--\eqref{assumptions-harn}, and for some point  $x_0 \in \Omega^+(u)\cup F(u)$,
\be\label{osc} (x_n+ a_0)^+ \leq u(x) \leq (x_n+ b_0)^+ \quad
\text{in $B_r(x_0) \subset \Omega,$}\ee
with 
$$b_0 - a_0 \leq \ep r, \qquad\ep \leq \bar \ep,$$ then
$$ (x_n+ a_1)^+ \leq u(x) \leq (x_n+ b_1)^+ \quad \text{in
$B_{r/40}(x_0)$},$$ with
$$a_0 \leq a_1 \leq b_1 \leq b_0, \quad
b_1 -  a_1\leq (1-c)\ep r, $$ and $0<c<1$ universal.
\end{thm}

\begin{proof} Assume without loss of generality that $x_0=0, r=1.$ 

 We call $q(x)=x_n+a_0$. Assumption \eqref{osc} gives that
\begin{equation}\label{trapped-by-q}
q^+(x)\leq u(x)\leq (q(x)+\varepsilon)^+\quad \:\:\mbox{in}\:\:B_1,
\end{equation}
since $b_0\le a_0+\varepsilon$. 
We distinguish three cases.

\vspace{1mm}

{\it Case 1.} $|a_0| < 1/20.$
We now distinguish two cases: $u(\hat{x}_0)\geq (q(\hat{x}_0)+\frac{\varepsilon}{2})^+$ or  $u(\hat{x}_0)\le (q(\hat{x}_0)+\frac{\varepsilon}{2})^+$, where $\hat{x}_0=\frac{1}{10}e_n$.

Assume that  
$$
u(\hat{x}_0)\geq (q(\hat{x}_0)+\frac{\varepsilon}{2})^+, \quad \hat{x}_0=\frac{1}{10}e_n,
$$
(the other case is treated similarly).
 Then, by Lemma \ref{pre-harnack}, if $\ep\le\bar{\ep}$,
\begin{equation*}
(q(x)+c\varepsilon)^+\leq u(x)\quad \:\:\mbox{in}\:\:\overline{B}_{\frac{1}{2}},
\end{equation*}
 for $0<c<1$ universal, which gives the desired improvement.

\vspace{1mm}

{\it Case 2.} $a_0 \le - 1/20.$ In this case it follows from \eqref{trapped-by-q} that, for $\ep<1/40$, $0$ belongs to the zero phase of
$(q(x)+\varepsilon)^+$, which implies that $0$ belongs to the zero phase of $u$. A contradiction.

\vspace{1mm}

{\it Case 3.} $a_0 \ge 1/20.$ In this case it follows from \eqref{osc} that 
$$B_{1/20} \subset B_1^+(u).$$
Then, denoting $\hat{u}=u-a_0$, we have
\begin{equation}\label{eq-p(x)hat}
\Delta_{p(x)}u = \Delta_{p(x)}\hat{u}=f\quad \text{ in } B_{1/20}.
\end{equation}
Observing that $||\hat{u}||_{L^{\infty}(B_1)}\le 2$ and recalling \eqref{assumptions-harn}, we obtain from the application of Theorem 1.1 
in \cite{Fan} to $\hat{u}$, that $u\in C^{1,\alpha}$ in $\overline{B}_{1/40},$  where $\alpha=\alpha (p_{\min},p_{\max},n)\in (0,1)$ and $||\nabla u||_{C^{\alpha}(\overline{B}_{1/40})}\leq C$, with $C=C (p_{\min},p_{\max},n)\geq 1.$

We now distinguish two cases: $u(0)-q(0)\ge \frac{\varepsilon}{2}$ or  $u(0)-q(0) \le \frac{\varepsilon}{2}$.

Assume that  
$$u(0)-q(0)\ge \frac{\varepsilon}{2},$$
(the other case is treated similarly). We will proceed as in the proof of Lemma \ref{pre-harnack}.

If $|\nabla u(0)|<\frac{1}{4}$, we argue as in  Case (i) of Lemma \ref{pre-harnack}, taking $\bar{x}_0=-r_2 e_n$. Here  $r_2>0$ is universal,
chosen as in that lemma, and such that we also have
$$B_{1/40}\subset\subset B_{r_4}(\bar{x}_0)\subset\subset B_{1/20},$$
for an appropriate chosen  universal $r_4>0$. We now take $r_3$ universal as in Lemma \ref{pre-harnack}, let
$$D:= B_{r_4}(\bar{x}_0)\setminus \overline{B_{r_3}(\bar{x}_0)},$$
and define $w$ in $D$ as in that lemma. Then, arguing as  in that proof,  we obtain
\begin{equation}\label{improv-case-3}
u(x)-q(x)\ge c_1{\varepsilon}\quad \text{ in } B_{1/40},
\end{equation}
with $0<c_1<1$, if $\ep\le\bar{\ep}$,  $\bar{\ep}$ and $c_1$ universal.

If $|\nabla u(0)|\ge \frac{1}{4}$, we proceed as in  Case (ii) of Lemma \ref{pre-harnack} and we consider the barrier $w$ in 
$$D:= B_{1/20}\setminus \overline{B_{r_0}},$$
with $r_0>0$ universal and small.  We obtain again \eqref{improv-case-3}, thus completing the proof.
\end{proof}

{}From Theorem $\ref{HI}$, with the same arguments employed in \cite{D}, we obtain the following estimate that will be crucial in the improvement of flatness procedure.

\begin{cor} \label{corollary}Let $u$ be as in Theorem \ref{HI}  satisfying \eqref{osc} for $r=1$. Then  in $B_1(x_0)$, $\tilde u_\ep(x) = \dfrac{u(x) - x_n }{\ep}$ 
 has a H\"older modulus of continuity at $x_0$, outside
the ball of radius $\ep/\bar \ep,$ i.e., for all $x \in \big(\Omega^+(u)\cup F(u)\big)\cap B_1(x_0)$, with $|x-x_0| \geq \ep/\bar\ep$,
$$|\tilde u_\ep(x) - \tilde u_\ep (x_0)| \leq C |x-x_0|^\gamma.
$$
Here $\bar \ep$ is as in Theorem \ref{HI}, and $C$ and $0<\gamma<1$ are universal.
\end{cor}

\section{Improvement of flatness} \label{sect-improv}

In this section we present the main improvement of flatness lemma.
Theorem \ref{flatmain1} will then be obtained by applying this lemma in an iterative way.
\begin{lem}[Improvement of flatness] \label{improv1}Let $u$  satisfy \eqref{fb} in $B_1$ and
\begin{equation}\label{ep-bound}
 \|f\|_{L^\infty(B_1)} \leq \ep^2, \quad ||g-1||_{L^{\infty}(B_1)}\le \ep^2, \quad  ||\nabla p||_{L^{\infty}(B_1)}\le {\ep}^{1+\theta},
\quad ||p-p_0||_{L^{\infty}(B_1)}\le \ep,
\end{equation}
for $0<\ep<1$, for some constant $0<\theta\le 1$. Suppose that
\begin{equation}\label{flat_1}(x_n -\ep)^+ \leq u(x) \leq (x_n +
\ep)^+ \quad \text{in $B_1,$} \quad 0\in F(u).
\end{equation}
If $0<r \leq r_0$ for $r_0$ universal, and $0<\ep \leq \ep_0$ for some $\ep_0$
depending on $r$, then
\begin{equation}\label{improvedflat_2_new}(x \cdot \nu -r \ep / 2)^+  \leq u(x) \leq
(x \cdot \nu +r \ep/2)^+ \quad \text{in $B_r,$}
\end{equation} with $|\nu|=1$ and $ |\nu - e_n| \leq \tilde C\ep$  for a
universal constant $\tilde C.$
\end{lem}
\begin{proof}We divide the proof of this lemma into 3 steps. We will use the following notation:
$$\Omega_{\rho}(u):=\big(B_1^+(u)\cup F(u)\big)\cap B_{\rho}.$$

 \vspace{2mm}

{\it Step 1: Compactness.} Fix $r \leq r_0$ with $r_0$ universal (the precise $r_0$ will be given in Step 3). Assume by contradiction that we
can find a sequence $\ep_k \rightarrow 0$ and a sequence $u_k$ of
solutions to \eqref{fb} in $B_1$ with  right hand side $f_k$, exponent $p_k$ and free boundary condition $g_k$ satisfying \eqref{ep-bound} with
$\ep=\ep_k$, such that $u_k$ satisfies \eqref{flat_1}, i.e.,
\begin{equation}\label{flat_k}(x_n -\ep_k)^+ \leq u_k(x) \leq (x_n + \ep_k)^+ \quad \text{for $x \in B_1$,  $0 \in F(u_k),$}
\end{equation} 
but $u_k$ does not satisfy the conclusion \eqref{improvedflat_2_new} of the lemma.

Set $$ \tilde{u}_{k}(x)= \dfrac{u_k(x) -  x_n}{\ep_k}, \quad x \in
\Omega_1(u_k).$$
Then, \eqref{flat_k} gives
\begin{equation}\label{flat_tilde**} -1 \leq \tilde{u}_{k}(x) \leq 1
\quad \text{for $x \in\Omega_1(u_k)$}.
\end{equation}

{}From Corollary \ref{corollary}, it follows that the function
$\tilde u_{k}$ satisfies \be\label{HC}|\tilde u_{k}(x) - \tilde
u_{k} (y)| \leq C |x-y|^\gamma,\ee for $C$ and $0<\gamma<1$ universal and
$$|x-y| \geq \ep_k/\bar\ep, \quad x,y \in\Omega_{1/2}(u_k).$$ From \eqref{flat_k} it clearly follows that
$F(u_k)$ converges to $B_1 \cap \{x_n=0\}$ in the Hausdorff
distance. This fact and \eqref{HC} together with Ascoli-Arzela
give that, as $\ep_k \rightarrow 0$, the graphs of the
$\tilde{u}_{k}$ over $\Omega_{1/2}(u_k)$ converge (up to a
subsequence) in the Hausdorff distance to the graph of a H\"older
continuous function $\tilde{u}$ over $B_{1/2}\cap\{x_n\ge 0\}$.

\vspace{2mm}

{\it Step 2: Limiting Solution.} We now show that $\tilde u$
solves the following linearized problem 
\begin{equation}\label{Neumann_p-in-1/2}
  \begin{cases}
    {\mathcal L}_{p_0} \tilde u=0 & \text{in $B_{1/2} \cap \{x_n > 0\}$}, \\
\tilde u_n=0 & \text{on $B_{1/2} \cap \{x_n =0\}$},
  \end{cases}\end{equation}
	in the sense of Definition \ref{def-neumann}. Here ${\mathcal L}_{p_0}$ is as in \eqref{Lp0}. 
	
	Let $P(x)$ be a quadratic polynomial  touching $\tilde u$  at $\bar{x}\in B_{1/2} \cap \{x_n \ge 0\}$ strictly from below. We need to show that
\begin{itemize}
\item[(i)] if $\bar{x}\in B_{1/2} \cap \{x_n > 0\}$ then ${\mathcal L}_{p_0}P\le 0$; 

\item[(ii)] if $\bar{x}\in B_{1/2} \cap \{x_n = 0\}$ then $P_n(\bar{x})\le 0$.
\end{itemize}

Since $\tilde{u}_{k}\to \tilde u$ in the sense specified above, there exist points $x_k\in \Omega_{1/2}(u_k)$, $x_k\to \bar{x}$ and constants
$c_k\to 0$ such that 
\begin{equation}\label{P-touches}
\tilde{u}_{k}(x_k)=P(x_k)+c_k
\end{equation}
and
\begin{equation}\label{P-below}
\tilde{u}_{k}\ge P+c_k\quad \text{ in  a neighborhood of } x_k.
\end{equation}
{}From the definition of $\tilde{u}_{k}$, \eqref{P-touches} and \eqref{P-below} read
\begin{equation*}
u_k(x_k)=Q_k(x_k)
\end{equation*}
and 
\begin{equation*}
{u}_{k}(x)\ge Q_k(x)\quad  \text{ in  a neighborhood of } x_k,
\end{equation*}
where
\begin{equation*}
Q_k(x)=\ep_k(P(x)+c_k)+x_n.
\end{equation*}
For notational simplicity we will drop the sub-index $k$ from $Q_k$.

We first notice that 
\begin{equation}\label{grad-Q}
\nabla Q=\ep_k \nabla P+e_n,
\end{equation}
thus, 
\begin{equation}\label{Qnot0}
\nabla Q(x_k) \neq 0,\quad \text{for $k$ large}.
\end{equation}

We now distinguish two cases.

(i) If $\bar{x}\in B_{1/2}\cap\{x_n>0\}$ then $x_k\in B^+_{1/2}(u_k)$ (for $k$ large).  Since $Q$ touches $u_k$ from below at $x_k$, and
$\nabla Q(x_k) \neq 0$, we get
\begin{equation*}
\begin{aligned}
\ep_k^2&\ge  f_k(x_k)\\
&\ge \Delta_{p_k(x_k)}Q(x_k)\\
&=|\nabla Q(x_k)|^{p_k(x_k)-2}\Delta Q + |\nabla
Q(x_k)|^{p_k(x_k)-4}(p_k(x_k)-2)\sum_{i,j=1}^n{Q}_{x_i}(x_k){Q}_{x_j}(x_k){Q}_{x_ix_j}\\
&\qquad+|\nabla Q(x_k)|^{p_k(x_k)-2}\langle \nabla p_k(x_k),\nabla Q(x_k)\rangle\log |\nabla Q(x_k)|\\
&=\ep_k|\nabla Q(x_k)|^{p_k(x_k)-2}\Delta P + \ep_k|\nabla
Q(x_k)|^{p_k(x_k)-4}(p_k(x_k)-2)\sum_{i,j=1}^n{Q}_{x_i}(x_k){Q}_{x_j}(x_k){P}_{x_ix_j}\\
&\qquad+|\nabla Q(x_k)|^{p_k(x_k)-2}\langle \nabla p_k(x_k),\nabla Q(x_k)\rangle\log |\nabla Q(x_k)|.
\end{aligned}
\end{equation*}
Using that $|\nabla p_k(x_k)|\le \ep_k$, we obtain
\begin{equation*}
\begin{aligned}
\ep_k&\ge |\nabla Q(x_k)|^{p_k(x_k)-2}\Delta P + |\nabla
Q(x_k)|^{p_k(x_k)-4}(p_k(x_k)-2)\sum_{i,j=1}^n{Q}_{x_i}(x_k){Q}_{x_j}(x_k){P}_{x_ix_j}\\
&\qquad-|\nabla Q(x_k)|^{p_k(x_k)-1}|\log |\nabla Q(x_k)||.
\end{aligned}
\end{equation*}
 Now, passing to the limit $k\to\infty$ and recalling that   
\begin{equation*}
\nabla Q(x_k)\to e_n,\qquad p_k(x_k)\to p_0, \qquad \ep_k\to 0,
\end{equation*}
we conclude that ${\mathcal L}_{p_0}P\le 0$ as desired.

(ii) If $\bar{x}\in B_{1/2}\cap\{x_n=0\}$, as observed in Remark  \ref{polyn-strict}, we can assume that ${\mathcal L}_{p_0}P> 0$. We claim that for $k$ large enough, $x_k\in F(u_k)$.  Otherwise $x_{k_j}\in B_{1/2}^+(u_{k_j})$ for a subsequence $k_j\to\infty$ and as in case (i), passing to the limit, we get
$${\mathcal L}_{p_0}P\le 0,$$
a contradiction. Thus, $x_k\in F(u_k)$ for $k$ large.

Since $Q^+$ touches $u_k$ from below at $x_k\in F(u_k)$ and \eqref{Qnot0} holds 
$$|\nabla Q(x_k)|\le g_k(x_k)\le 1+\ep_k^2,$$
which, by \eqref{grad-Q}, gives
$$|\nabla Q(x_k)|^2 = \ep_k^2|\nabla P(x_k)|^2+1+2\ep_kP_n(x_k)\le 1+3\ep_k^2.$$
Thus, after division by $\ep_k$,
$$ \ep_k|\nabla P(x_k)|^2-3\ep_k+2P_n(x_k)\le 0.$$
Passing to the limit as $k\to\infty$, we obtain $P_n(\bar{x})\le 0$ as desired.

\vspace{2mm}

{\it Step 3: Improvement of flatness.} From the previous step, $\tilde u$ solves \eqref{Neumann_p-in-1/2} and from \eqref{flat_tilde**},
\begin{equation*}-1 \leq \tilde{u}(x) \leq 1
\quad \text{in $B_{1/2}\cap\{x_n\ge 0\}$}.
\end{equation*}

{}From Theorem \ref{linearreg} and the bound above we find that, for the given $r$,
$$|\tilde u(x) - \tilde u(0) - \nabla \tilde u(0)\cdot x| \leq C_0 r^2 \quad \text{in } B_r\cap\{x_n\ge 0\},$$
if $r_0\le 1/4$, for a universal constant $C_0$. In particular, since  $\tilde u(0)=0$ and also $\tilde u_n(0)=0$, we obtain 
$$ x'\cdot \tilde \nu - C_0 r^2\le \tilde u(x)\le x'\cdot \tilde \nu + C_0 r^2 \quad \text{in } B_r\cap\{x_n\ge 0\},$$
where $x'=(x_1,\cdots,x_{n-1})$, $\tilde\nu=\nabla_{x'}\tilde u(0)$ and $|\tilde \nu|\le C_0$. Therefore, for $k$ large enough we get
$$ x'\cdot \tilde \nu - C_1 r^2\le \tilde u_k(x)\le x'\cdot \tilde \nu + C_1 r^2 \quad \text{in } \Omega_r(u_k),$$
for a universal constant $C_1$. From the definition of $\tilde u_k$ the inequality above reads
\begin{equation}\label{improve-u_k}
\ep_k x'\cdot \tilde \nu +x_n - \ep_k C_1 r^2\le u_k(x)\le \ep_k x'\cdot \tilde \nu +x_n + \ep_k C_1 r^2 \quad \text{in } \Omega_r(u_k).
\end{equation}
We next set
$$ \nu_k = \frac{1}{\sqrt {1+ \ep_k^2 |\tilde\nu|^2}}(e_n + \ep_k (\tilde\nu, 0)).$$
Then,
$$|\nu_k|=1,\qquad |\nu_k-e_n|\le \tilde C \ep_k,$$
and 
$$\nu_k =e_n + \ep_k (\tilde\nu, 0) + \ep_k^2 \tau, \quad |\tau|\leq \tilde C,$$
with $\tilde C$ universal. We now deduce from \eqref{improve-u_k}
\begin{equation*}
x\cdot \nu_k -\ep^2_k \tilde C r - \ep_k C_1 r^2\le u_k(x)\le x\cdot \nu_k +\ep^2_k \tilde C r + \ep_k C_1 r^2 \quad \text{in } \Omega_r(u_k).
\end{equation*}
If we fix $r_0$ satisfying $C_1 r_0\le 1/4$ and we take $k$ large enough so that $\ep_k \tilde C\le 1/4$ , we get
\begin{equation*}
x\cdot \nu_k -\ep_k  r/2 \le u_k(x)\le x\cdot \nu_k +\ep_k r/2 \quad \text{in } \Omega_r(u_k).
\end{equation*}
Recalling \eqref{flat_k}, we obtain for large $k$
\begin{equation*}
(x\cdot \nu_k -\ep_k  r/2)^+ \le u_k(x)\le (x\cdot \nu_k +\ep_k r/2)^+ \quad \text{in } B_r,
\end{equation*}
thus $u_k$ satisfies the conclusion \eqref{improvedflat_2_new} of the lemma, a contradiction.
\end{proof}

\section{Regularity of the free boundary}\label{section7}

In this section we finally prove our main result, namely, Theorem \ref{flatmain1}.

\begin{proof}[\bf Proof of Theorem \ref{flatmain1}]
Let $u$ be a viscosity solution to \eqref{fb}
in $B_1$ with $0\in F(u),$ $g(0)=1$ and $p(0)=p_0.$   Consider the sequence 
$$u_k(x) = \frac{1}{\rho_k} u(\rho_k x), \quad x\in B_1, $$
with $\rho_k = \bar r^k$, $k=0,1,\cdots$, for a fixed $\bar r$ such that 
$$\bar r^{\beta} \le 1/4,\quad \bar r\le r_0,$$
with $r_0$ the universal constant in Lemma \ref{improv1}, taking $\theta=1$ in \eqref{ep-bound}.

Each $u_k$ is a solution to \eqref{fb} with right hand side $f_k(x) = \rho_k f(\rho_k x)$, exponent $p_k(x) =  p(\rho_k x)$, and 
free boundary condition  $g_k(x) =  g(\rho_k x)$. For the chosen $\bar r$, by taking $\bar\ep=\ep_0(\bar r)^2$, the assumption \eqref{ep-bound} holds for $\ep=\ep_k=2^{-k}\ep_0(\bar r)$. Indeed, in $B_1$, in view of \eqref{pflat},
\begin{equation*}
\begin{aligned}
|f_k(x)|&\le ||f||_{\infty} \,\rho_k\le \bar\ep {\bar r}^k\le \ep_k^2,\\
|g_k(x)-1|&=|g(\rho_k x)-g(0)|\le [g]_{0,\beta}\, {\rho_k}^{\beta}\le\bar\ep {\bar r}^{k\beta}\le \ep_k^2,\\
|\nabla p_k(x)|&\le ||\nabla p||_{\infty}\,\rho_k\le \bar\ep {\bar r}^k\le \ep_k^2,\\
|p_k(x)-p_0|&=|p(\rho_k x)-p(0)|\le ||\nabla p||_{\infty}\,\rho_k\le \bar\ep {\bar r}^k\le \ep_k^2.
\end{aligned}
\end{equation*}

The hypothesis \eqref{cflat} guarantees that for $k=0$ also the flatness assumption \eqref{flat_1} in Lemma \ref{improv1} is satisfied by
$u_0$. Then it easily follows, by applying inductively Lemma \ref{improv1}, that each $u_k$ is $\ep_k-$flat in $B_1$ in the sense of 
\eqref{flat_1}, in the direction $\nu_k$, with $|\nu_k|=1,$ $ |\nu_k - \nu_{k+1}| \leq \tilde C \ep_k$ ($\nu_0=e_n$). Now, a standard iteration argument gives the desired statement.
\end{proof}

\medskip

\appendix

\section{Lebesgue and Sobolev spaces with variable
exponent} \label{appA1}

Let $p :\Omega \to  [1,\infty)$ be a measurable bounded function,
called a variable exponent on $\Omega$, and denote $p_{\max} = {\rm
ess sup} \,p(x)$ and $p_{\min} = {\rm ess inf} \,p(x)$. The variable exponent Lebesgue space $L^{p(\cdot)}(\Omega)$ is defined as the set of all measurable functions $u :\Omega \to \R$ for which
the modular $\varrho_{p(\cdot)}(u) = \int_{\Omega} |u(x)|^{p(x)}\,
dx$ is finite. The Luxemburg norm on this space is defined by
$$
\|u\|_{L^{p(\cdot)}(\Omega)} = \|u\|_{p(\cdot)}  = \inf\{\lambda >
0: \varrho_{p(\cdot)}(u/\lambda)\leq 1 \}.
$$

This norm makes $L^{p(\cdot)}(\Omega)$ a Banach space.

There holds the following relation between $\varrho_{p(\cdot)}(u)$
and $\|u\|_{L^{p(\cdot)}}$:
\begin{align*}
\min\Big\{\Big(\int_{\Omega} |u|^{p(x)}\, dx\Big)
^{1/{p_{\min}}},& \Big(\int_{\Omega} |u|^{p(x)}\, dx\Big)
^{1/{p_{\max}}}\Big\}\le\|u\|_{L^{p(\cdot)}(\Omega)}\\
 &\leq  \max\Big\{\Big(\int_{\Omega} |u|^{p(x)}\, dx\Big)
^{1/{p_{\min}}}, \Big(\int_{\Omega} |u|^{p(x)}\, dx\Big)
^{1/{p_{\max}}}\Big\}.
\end{align*}

Moreover, the dual of $L^{p(\cdot)}(\Omega)$ is
$L^{p'(\cdot)}(\Omega)$ with $\frac{1}{p(x)}+\frac{1}{p'(x)}=1$.

 $W^{1,p(\cdot)}(\Omega)$ denotes the space of measurable
functions $u$ such that $u$ and the distributional derivative
$\nabla u$ are in $L^{p(\cdot)}(\Omega)$. The norm

$$
\|u\|_{1,p(\cdot)}:= \|u\|_{p(\cdot)} + \| |\nabla u|
\|_{p(\cdot)}
$$
makes $W^{1,p(\cdot)}(\Omega)$ a Banach space.

The space $W_0^{1,p(\cdot)}(\Omega)$ is defined as the closure of
the $C_0^{\infty}(\Omega)$ in $W^{1,p(\cdot)}(\Omega)$.

For  further details on these spaces, see \cite{DHHR},
\cite{KR}, \cite{RaRe} and their references.

\section*{Acknowledgment }
The authors wish to thank Sandro Salsa for very interesting discussions about the subject of this paper.


\begin{thebibliography}{9999}

\bibitem[AMS]{AMS} R. Aboulaich, D. Meskine, A. Souissi, {\it New diffusion models in image processing}, Comput.
Math. Appl. 56 (2008) 874--882.
\bibitem[AM]{AM} E. Acerbi, G. Mingione, {\it Regularity results for a class of functionals with non-standard growth,} Arch. Ration. Mech. Anal. 156 (2) (2001) 121--140. 

\bibitem[AC]{AC} H. W. Alt,
L. A. Caffarelli, {\it Existence and regularity for a minimum problem
with free boundary}, J. Reine Angew. Math 325
(1981) 105--144.


\bibitem[ACF]{ACF}
H.~W. Alt, L.~A. Caffarelli, A.~Friedman, {\it A free boundary
problem for
  quasilinear elliptic equations}, Ann. Sc. Norm. Super. Pisa Cl. Sci. (4)
  11 (1) (1984)  1--44.


\bibitem[AR]{AR} S. N. Antontsev, J. F. Rodrigues, {\it On stationary thermo-rheological viscous flows}, Ann. Univ. Ferrara, Sez. VII, Sci. Mat. 
 52 (1) (2006) 19--36.
\bibitem[ART]{ART}  D. J. Ara\'ujo, G. C. Ricarte, E. V. Teixeira, {\it Singularly perturbed equations of degenerate type,} Ann. Inst. H. Poincar\'e Anal. Non Lin\'eaire 34 (3) (2017)  655--678.

\bibitem[AF]{AF} R. Argiolas, F. Ferrari, {\it Flat free boundaries regularity in two-phase problems for a class of fully nonlinear elliptic operators with variable coefficients,} Interfaces Free Bound. 11 (2) (2009) 177-199.

\bibitem[C1]{C1} L. A. Caffarelli, {\it A Harnack
inequality approach to the regularity of free boundaries. Part I:
Lipschitz free boundaries are $C^{1,\alpha}$}, Rev. Mat.
Iberoamericana 3 (2) (1987) 139--162.
\bibitem[C2]{C2} L. A. Caffarelli, {\it A Harnack
inequality approach to the regularity of free boundaries. Part II:
Flat free boundaries are Lipschitz}, Comm. Pure Appl. Math.
42 (1) (1989) 55--78.



\bibitem[CFS]{CFS} M. C.  Cerutti, F. Ferrari, S. Salsa, {\it Two phase
problems for linear elliptic operators with variable coefficients:
Lipschitz free boundaries are $C^{1,\gamma}$}, Archive for
Rational Mechanics and Analysis 171 (3) (2004)  329 - 348.



\bibitem[CL]{CL} S. Challal, A. Lyaghfouri, {\it Second order regularity for the $p(x)$-Laplace operator}, Math. Nachr. 284 (10)
 (2011) 1270--1279.

\bibitem[CLR]{CLR} Y. Chen, S. Levine, M. Rao, {\it Variable exponent, linear growth functionals in image restoration},
SIAM J. Appl. Math. 66 (2006) 1383--1406.


\bibitem[CIL]{CIL} M. G. Crandall, H. Ishii, P. L. Lions, {\it User’s guide to viscosity solutions of second order partial differential equations}, Bull. Amer. Math. Soc. (N.S.) 27 (1) (1992)  1--67.

\bibitem[DP]{DP}
D.~Danielli, A.~Petrosyan, \emph{A minimum problem with free
boundary for a
  degenerate quasilinear operator}, Calc. Var. Partial Differential Equations
  23 (1) (2005)  97--124.


\bibitem[D]{D} D. De Silva, {\it Free boundary regularity for a problem with right hand side}, Interfaces and free boundaries 13 (2011) 223--238.
\bibitem[DFS1]{DFS1} D. De Silva, F. Ferrari, S. Salsa, {\it Two-phase problems with distributed sources: regularity of the free boundary.} Anal. PDE 7 (2) (2014) 267--310. 

\bibitem[DFS2]{DFS2} D. De Silva, F. Ferrari, S. Salsa, {\it Free boundary regularity for fully nonlinear non-homogeneous two-phase problems.} J. Math. Pures Appl. (9) 103 (3) (2015) 658--694.

\bibitem[DFS3]{DFS3} D. De Silva, F. Ferrari, S. Salsa, {\it Regularity of higher order in two-phase free boundary problems.} Trans. Amer. Math. Soc. 371 (5) (2019) 3691--3720.

\bibitem[DHHR]{DHHR} L. Diening, P. Harjulehto, P. Hasto, M. Ruzicka, {\it Lebesgue and Sobolev Spaces with variable exponents}, Lecture Notes in Mathematics 2017, Springer,  2011.


\bibitem[Fa]{Fan} X. Fan, {\it Global $C^{1,\alpha}$ regularity for variable exponent elliptic equations in divergence form}, J. Differential Equations 235 (2007) 397--417.


\bibitem[FZ]{FanZ} X. Fan, D. Zhao, {\it A class of De Giorgi type and H\"older continuity}, Nonlinear Analysis 36 (1999)  295--318.



\bibitem[F1]{F1} M. Feldman, {\it Regularity for nonisotropic two-phase problems with Lipschitz free boundaries}, Differential Integral Equations 10  (6) (1997) 1171--1179.

\bibitem[F2]{F2} M. Feldman, {\it Regularity of Lipschitz free boundaries in two-phase problems for fully nonlinear elliptic equations}, Indiana Univ. Math. J. 50 (3) (2001) 1171--1200.

\bibitem[FMW]{FMW} J. Fernandez  Bonder, S. Mart\'{\i}nez, N. Wolanski, {\it A free boundary problem for the $p(x)$-Laplacian}, Nonlinear Anal.
72 (2010) 1078--1103.


\bibitem[Fe1]{Fe1} F. Ferrari, {\it Two-phase problems for a class of fully nonlinear elliptic operators, Lipschitz free boundaries are $C^{1,\gamma}$},
Amer. J. Math. 128 (2006) 541--571.




\bibitem[FL]{FL}  F. Ferrari, C. Lederman, {\it Regularity of Lipschitz free boundaries for a $p(x)$-Laplacian problem with right hand side}, preprint.

 
\bibitem[FS1]{FS1} F. Ferrari, S. Salsa, {\it Regularity of the free boundary in two-phase problems for elliptic operators},
 Adv. Math. 214 (2007) 288--322.

\bibitem[FS2]{FS2} F. Ferrari, S. Salsa, {\it Subsolutions of elliptic operators in divergence form and application to two-phase free boundary problems}, Bound. Value Probl. 2007, art. ID 57049, 21pp.

\bibitem[GT]{GT}
D. Gilbarg, N. S. Trudinger, \emph{Elliptic Partial Differential
Equations
  of Second Order}, Grundlehren der Mathematischen Wissenschaften [Fundamental
  Principles of Mathematical Sciences], vol. 224, Springer-Verlag, Berlin,
  1983 (3rd edition).


  \bibitem[GS]{GS} B. Gustafsson, H. Shahgholian, {\it Existence and
geometric properties of solutions of a free boundary problem in
potential theory}, J. Reine Angew. Math. 473 (1996)
137--179.


\bibitem[JJ]{JJ} V. Julin, P. Juutinen, {\it A new proof for the equivalence of weak and viscosity solutions for the $p$-Laplace equation}, Communications in PDE 37 (5) (2012) 934 -- 946.

\bibitem[JLM]{JLM} P. Juutinen, P. Lindqvist, J. Manfredi, {\it On the equivalence of viscosity solutions and weak solutions for a quasi-linear equation}. SIAM J. Math. Anal. 33 (3) (2001) 699--717.


 \bibitem[JLP]{JLP} P. Juutinen, T. Lukkari, M. Parviainen, {\it Equivalence of viscosity and weak solutions for the
$p(x)$-Laplacian}. Ann. Inst. H. Poincare Anal. Non Lineaire 27 (6) (2010) 1471--1487.

\bibitem[KR]{KR}
O. Kov\'a\v{c}ik, J. R\'akosn{\'i}k, \emph{On spaces ${L}^{p(x)}$ and
  ${W}^{k,p(x)}$}, Czechoslovak Math. J  41 (1991) 592--618.

\bibitem[Le]{Le}
C.~Lederman, \emph{A free boundary problem with a volume penalization}, Ann.
  Sc. Norm. Super. Pisa Cl. Sci. (4) 23 (2) (1996)  249--300.


\bibitem[LW1]{LW1} C.  Lederman, N. Wolanski, {\it An inhomogeneous singular perturbation problem for the $p(x)$-Laplacian,} Nonlinear Anal. 138 (2016) 300--325.
\bibitem[LW2]{LW2} C.  Lederman, N. Wolanski, {\it Weak solutions and regularity of the interface in an inhomogeneous free boundary problem for the $p(x)$-Laplacian,} Interfaces Free Bound. 19 (2) (2017) 201--241.
\bibitem[LW3]{LW3} C. Lederman, N. Wolanski, {\it Inhomogeneous minimization problems for the $p(x)$-Laplacian,} J. Math. Anal. Appl. 475 (1) (2019) 423--463.
\bibitem[LW4]{LW4} C. Lederman, N. Wolanski, {\it An optimization problem with volume constraint for an inhomogeneous operator with nonstandard growth},  Discrete Contin. Dyn. Syst.  41 (6) (2021) 2907--2946. 
\bibitem[LR]{LR} R. Leit$\tilde{a}$o, G. Ricarte, 
 {\it Free boundary regularity for a degenerate problem with right hand side},
Interfaces Free Bound. 20 (2018) 577--595.
\bibitem[LT]{LT} R. Leit$\tilde{a}$o, E. V.  Teixeira, {\it Regularity and geometric estimates for minima of discontinuous functionals}, Rev. Mat. Iberoam. 31 (1) (2015) 69--108. 

\bibitem[LN1]{LN1} J. Lewis, K. Nystr{\"o}m, {\it Regularity of Lipschitz free boundaries in two phase problems for the $p$-Laplace operator,} Adv. in Math. 225 (2010) 2565-2597. 
\bibitem[LN2]{LN2} J. Lewis, K. Nystr{\"o}m K., {\it Regularity of flat free boundaries in two-phase problems for the $p$-Laplace operator,} Ann. Inst. H. Poincar\'e Anal. Non Lin\'aire 29 (1) (2012) 83--108.

\bibitem[MW]{MW}
S.~Mart\'{\i}nez, N.~Wolanski, \emph{A minimum problem with free boundary in
  {O}rlicz spaces}, Adv. Math. 218 (6) (2008)
  1914--1971.


\bibitem[MO]{MO} M. Medina, P. Ochoa, {\it On the viscosity and weak solutions for non-homogeneous $p$-Laplace equations}. Adv. in Nonlinear Anal. 8 (1) (2019) 468--481.


\bibitem[MS]{MS} E. Milakis, L. Silvestre, {\it  Regularity for fully nonlinear elliptic equations with Neumann boundary data}, Comm. in Partial Diff. Equations 31 (2006)  1227--1252.

\bibitem[RR]{RaRe} V. D. Radulescu, D. D. Repovs, {\it Partial differential equations with variable exponents: variational methods and qualitative analysis}, Monographs and Research Notes in Mathematics, Book 9. Chapman \& Hall / CRC Press, Boca Raton, FL, 2015.


\bibitem[RT]{RT}  G. C. Ricarte, E. V. Teixeira, {\it Fully nonlinear singularly perturbed equations and asymptotic free boundaries}, J. Funct. Anal. 261 (2011) 1624--1673.


\bibitem[R]{R} M. Ruzicka, {\it Electrorheological Fluids: Modeling and
Mathematical Theory}, Springer-Verlag, Berlin, 2000.


\bibitem[S]{S}  O. Savin, {\it Small perturbation solutions for elliptic equations}. Comm. Partial Differential Equations 32 (4-6) (2007) 557--578.

\bibitem[W1]{W1}  P. Y. Wang, {\it Regularity of free boundaries of two-phase problems for fully nonlinear elliptic equations of second order. I. Lipschitz free boundaries are $C^{1,\alpha}$}, Comm. Pure Appl. Math. 53 (2000) 799--810.
\bibitem[W2]{W2}  P. Y. Wang, {\it Regularity of free boundaries of two-phase problems for fully nonlinear elliptic equations of second order. II. Flat free boundaries are Lipschitz}, Comm. Partial Differential Equations 27 (2002) 1497--1514.
\bibitem[Wo]{Wo} N. Wolanski, {\it Local bounds, Harnack inequality and H\"older continuity for divergence type elliptic equations with non-standard growth},
Rev. Un. Mat. Argentina 56 (1) (2015) 73--105.

\bibitem[Z1]{Z1} V. V. Zhikov, {\it Averaging of functionals of the calculus of variations and elasticity theory}, Math. USSR. Izv. 29 (1)  (1987) 33--66.


\bibitem[Z2]{Z2} V. V. Zhikov, {\it Solvability of the three-dimensional thermistor problem}, Tr. Mat. Inst. Steklova D (Differ. Uravn. i Din. Sist.) 261 (2008) 101–114.


\end{thebibliography}
\end{document}